\documentclass[12pt]{amsart}
\usepackage{amsmath,amstext,amsbsy,amssymb,amscd}
\usepackage{amsmath}
\usepackage{amsxtra}
\usepackage{amscd}
\usepackage{amsthm}
\usepackage{amsfonts}
\usepackage{amssymb}
\usepackage{eucal}
\usepackage{color}
\usepackage[enableskew,vcentermath]{youngtab}

\newtheorem{thm}{Theorem}[section]
\newtheorem{theorem}{Theorem}[section]
\theoremstyle{plain}

\newtheorem{proposition}[thm]{Proposition}

\theoremstyle{definition}

\newtheorem{remark}[thm]{Remark}

\definecolor{A}{rgb}{.75,1,.75}

\newcommand{\rarrow}{\rightarrow}
\newcommand{\cd}{{\, \cdot \, }}

\newcommand{\calK}{{\mathcal{K}}}

\newcommand{\dis}{\displaystyle}
\newcommand{\disp}{\displaystyle}
\newcommand{\q}{\bf q}
\newcommand{\p}{\bf p}

\newcommand{\f}{\bf f}
\newcommand{\be}{\begin{equation}}
\newcommand{\ee}{\end{equation}}
\newcommand{\parens}[1]{\left(#1\strut\right)}

\newcommand{\pa}{\partial}

{\vskip-\lastskip\medskip
  \noindent
  {\em #1.}\enspace
  }%
{\qed\par\medskip
  }

\title[Recover All Coefficients]{Recover All Coefficients in Second-order Hyperbolic Equations from Finite Sets of Boundary Measurements}
\author[Shitao Liu, Antonio Pierrottet and Scott Scruggs]{Shitao Liu, Antonio Pierrottet and Scott Scruggs}
\address{School of Mathematical and Statistical Sciences, Clemson University, Clemson, SC 29634}
\email{liul@clemson.edu}
\email{srscrug@g.clemson.edu}
\email{ampierr@g.clemson.edu}

\date{\today}

\begin{document}

\maketitle
\begin{abstract}
We consider the inverse hyperbolic problem of recovering all spatial dependent coefficients, which are the wave speed, the damping coefficient, potential coefficient and gradient coefficient, in a second-order hyperbolic equation defined on an open bounded domain with smooth enough boundary. We show that by appropriately selecting finite pairs of initial conditions we can uniquely and Lipschitz stably recover all those coefficients from the corresponding boundary measurements of their solutions. The proofs are based on sharp Carleman estimate, continuous observability inequality and regularity theory for general second-order hyperbolic equations. \end{abstract}




\section{Introduction and Main Results}\label{sec1}

Let $\Omega\subset\mathbb{R}^n$, $n\geq 2$, be an open bounded domain with smooth enough (e.g., $C^2$) boundary $\Gamma=\pa\Omega=\overline{\Gamma_0\cup\Gamma_1}$,
where $\Gamma_0\cap\Gamma_1=\emptyset$. We refer $\Gamma_1$ as the observed part of the boundary where the measurements are taken, and $\Gamma_0$ as the unobserved part of the boundary.
We consider the following general second-order hyperbolic equation for $w=w(x,t)$ defined on $Q=\Omega\times[-T,T]$, along with initial conditions $\{w_0, w_1\}$ and Dirichlet boundary condition $h$ on $\Sigma=\Gamma\times [-T,T]$ that are given in appropriate function spaces:
\begin{equation}\label{1}
\begin{cases}
 w_{tt} - c^2(x)\Delta w +q_1(x)w_t+q_0(x)w+ {\bf q}(x)\cdot\nabla w = 0 & \mbox{in } Q \\[2mm]
 w\left(x,0\right) = w_0(x); \ w_t\left(x,0\right) = w_1(x) & \mbox{in } \Omega \\[2mm]
 w(x,t) = h(x,t) & \mbox{in } \Sigma.
\end{cases}
\end{equation}
Here the wave speed $c(x)$ satisfies
\begin{equation*}
c\in\mathcal{C}=\{c\in C^1(\Omega): c_0^{-1}\leq c(x)\leq c_0, \ \mbox{for some} \ c_0>0\}
\end{equation*}
$q_1\in L^\infty(\Omega)$, $q_0\in L^\infty(\Omega)$, and ${\bf q}\in \left(L^\infty(\Omega)\right)^n$ are the damping, potential, and gradient coefficients, respectively.  

We then consider the following inverse problem for the system (\ref{1}): Recover all together the wave speed $c(x)$, the damping coefficient $q_1(x)$, the potential coefficient $q_0(x)$, and the gradient coefficient ${\bf q}(x)$ from measurements of Neumann boundary traces of the solution $w = w(w_0,w_1,h,c,q_1,q_0,{\bf q})$ over the observed part $\Gamma_1$ of the boundary and over the time interval $[-T,T]$. Of course here $T>0$ should be sufficiently large due to the finite propagation speed of the system \eqref{1}. In addition, to make the observed part $\Gamma_1$ of the boundary more precise, in this paper we assume the following standard geometrical assumptions on the domain $\Omega$ and the unobserved part of the boundary $\Gamma_0$:

\smallskip
(A.1) There exists a strictly convex function $d: \overline{\Omega} \rarrow \mathbb{R}$ in the metric $g=c^{-2}(x)dx^2$, and of class $C^3(\overline{\Omega})$, such that the following two properties hold true (through translation and rescaling if necessary):

\smallskip
(i) The normal derivative of $d$ on the unobserved part $\Gamma_0$ of the boundary is non-positive. Namely,
\begin{equation*}
\frac{\pa d}{\pa\nu}=\langle Dd(x), \nu(x)\rangle\leq 0, \quad  \forall x\in\Gamma_0,
\end{equation*}
where $Dd=\nabla_g d$ is the gradient vector field on $\Omega$ with respect to $g$.

\medskip
(ii) \begin{equation*} D^2d(X,X) = \langle D_X(Dd), X\rangle_g\geq 2|X|_g^2, \ \forall X\in M_x, \ \min_{x\in\overline{\Omega}}d(x)=m_0>0\end{equation*}
where $D^2d$ is the Hessian of $d$ (a second-order tensor) and $M_x$ is the tangent space at $x\in\Omega$.

\smallskip
(A.2) $d(x)$ has no critical point on $\overline{\Omega}$. In other words,
\begin{equation*}
\inf_{x\in\overline{\Omega}}|Dd|>0, \ \mbox{so that we may take} \ \inf_{x\in\overline{\Omega}}\frac{|Dd|^2}{d}>4.
\end{equation*}

\begin{remark}\label{rem2}
The geometrical assumptions above permit the construction of a vector field that enables a pseudo-convex function necessary for allowing a Carleman estimate containing no lower-order terms for the general second-order equation (\ref{1}) (see Section 2).  These assumptions are first formulated in \cite{LTZ2000} under the framework of a Euclidean metric, with \cite{TY2002} employing them under the more general Riemannian framework.  For examples and detailed illustrations of large general classes of domains $\{\Omega, \Gamma_1, \Gamma_0\}$ satisfying the aforementioned assumptions we refer to \cite[Appendix B]{TY2002}. One canonical example is to take $d(x)=|x-x_0|^2$, with $x_0$ being a point outside $\overline{\Omega}$, if the wave speed $c$ satisfies $\left|\frac{\nabla c(x)\cdot(x-x_0)}{2c(x)}\right|\leq r_c<1$ for some $r_c\in(0,1)$.
\end{remark}

The classical inverse hyperbolic problems usually involve recovering a single unknown coefficient, typically the damping coefficient or the potential coefficient, from a {\it single} boundary measurement of the solution. To some extent, those setup are expected since the unknown coefficient, whether it is the damping or the potential one, depends on $n$ independent variables and the corresponding boundary measurement also depends on $n$ free variables. In fact, under proper conditions it is even possible to recover both potential and damping coefficients in one shot by just one single boundary measurement \cite{LiuT2011}. 
In the case of a gradient coefficient, the unknown function is vector-valued and containing $n$ different real-valued functions. Hence a single measurement does not seem to be sufficient to recover all of them, which is probably why such problem is much less studied in the literature. Nevertheless, it is possible to recover the coefficient by properly making $n$ sets of boundary measurements \cite{Jellali2006}. Last, in the case of recovering the variable unknown wave speed, since the unknown function is at the principle order level, one typically needs to rewrite the hyperbolic equation as a Riemannian wave equation so that the principle part becomes constant coefficients on an appropriate Riemannian manifold \cite{B2004, Liu2013}.  

In this paper, we seek to recover all together the aformentioned coefficients in the second-order hyperbolic equation (\ref{1}). To the best of our knowledge, this is the first paper that addresses the uniqueness and stability of recovering all these coefficients at once through finitely many boundary measurements. Note that all together these coefficients contain a total of $n+3$ unknown functions, so naturally one may expect to be able to recover them by making $n+3$ sets of boundary measurements. This is entirely possible to do following the ideas in this paper (see Remark (1) in Section 4). Nevertheless, in the following we will show that by appropriately choosing $\lfloor{\frac{n+4}{2}}\rfloor$\footnote{Here $\lfloor\cdot\rfloor$ denotes the usual floor function.} pairs of initial conditions $\{w_0,w_1\}$ and a boundary condition $h$, we can uniquely and Lipschitz stably recover the coefficients $c$, $q_1,$ $q_0,$ and $\q$ all at once from the corresponding Neumann boundary measurements of their solutions. The precise results are stated in Theorem \ref{Th1} and Theorem \ref{Th2} below.

As mentioned above recovering a single coefficient from a single boundary measurement is a standard formulation in inverse hyperbolic problems and such problem has been studied extensively in the literature. 
Here we only mention the monographs and lecture notes \cite{B-Y2017,Isakov2000,Isakov2006,Klibanov2004,Klibanov2021,LRS1986,Liu-T2013} and refer to the substantial lists of references therein. The standard approach for this type of inverse hyperbolic problems typically involves using Carleman-type estimates for the second-order hyperbolic equations. To certain extent, such methods can all be seen as variations or improvements of the so called Bukhgeim--Klibanov (BK) method which was originated in the seminal paper \cite{BK1981}. Our approach to solve the present inverse problem also relies on a sharp Carleman estimate for general second-order hyperbolic equations and in particular a post Carleman estimate route that was introduced by Isakov in \cite[Theorem~8.2.2]{Isakov2006}. Another standard feature of the BK method is the need of certain positivity assumptions on the initial conditions. Of course the precise assumption depends on what coefficient(s) one is trying to recover.  

On the other hand, let us also mention that there is another standard formulation of inverse hyperbolic problems that usually does not require positivity on the initial conditions. In this formulation, one tries to recover information of second-order hyperbolic equations from all possible boundary measurements, which are often modeled by the Dirichlet to Neumann or Neumann to Dirichlet operator. In particular, in this case it is possible to recover all coefficients in system (\ref{1}) up to natural gauge transformations \cite{KL2002}, following the powerful Boundary Control (BC) method developed by Belishev \cite{Be1987}. For more inverse hyperbolic problems with infinitely many measurements and the BC method, we refer to the review paper \cite{Be1997} and the monograph \cite{KKL2000}.

\medskip
Let us now state the main theorems in this paper. 

\begin{theorem}\label{Th1}
Under the geometrical assumptions (A.1) and (A.2) and let 
\begin{equation}\label{largetime}
T > T_0 = 2\dis\sqrt{\max_{x\in\overline{\Omega}}d(x)}.
\end{equation} 
Suppose the initial and boundary conditions are in the following function spaces 
\be\label{regularity}
\{w_0, w_1\}\in H^{\gamma+1}(\Omega)\times H^{\gamma}(\Omega), \ h\in H^{{\gamma}+1}(\Sigma), \ \mbox{where} \ \gamma>\frac{n}{2}+4
\ee
along with all compatibility conditions (trace coincidence) which make sense. In addition, depending on the dimension $n$ of the space, we assume the following positivity condition: There exists $r_0>0$ such that 

\medskip
Case I: If $n$ is odd, i.e., $n=2m+1$ for some $m\in\mathbb{N}$, then we choose $m+2$ pairs of initial conditions $\{w_0^{(i)}, w_1^{(i)}\}$, $i=1, \dots, m+2$, and a boundary condition $h$ so that they satisfy (\ref{regularity}) and
\begin{equation}\label{positivity1'}
|\det{W(x)}|\geq r_0, \ a.e.\  x\in\Omega
\end{equation}
where $W(x)$ is the $(n+3)\times(n+3)$ matrix defined by

\begin{equation}\label{8}
W(x) =
\begin{bmatrix}
w_0^{(1)}(x) & w_1^{(1)}(x)  & \dis\pa_{x_1} w_0^{(1)}(x)  & \cdots & \dis\pa_{x_n} w_0^{(1)}(x)  &\Delta w_0^{(1)}(x)  \\[2mm]
w_1^{(1)}(x)  & w_{tt}^{(1)}(x)  &\pa_{x_1}w_1^{(1)}(x) & \cdots &\pa_{x_n}w_1^{(1)}(x)  &\Delta w_1^{(1)}(x)  \\[2mm]
\vdots & \vdots & \vdots & \ddots & \vdots & \vdots \\
w_0^{(m+2)}(x)  & w_1^{(m+2)}(x)  &\dis\pa_{x_1} w_0^{(m+2)}(x)  & \cdots & \dis\pa_{x_n} w_0^{(m+2)}(x)  & \Delta w_0^{(m+2)}(x) \\[2mm]
w_1^{(m+2)}(x)  & w_{tt}^{(m+2)}(x)  &\pa_{x_1}w_1^{(m+2)}(x) & \cdots &\pa_{x_n}w_1^{(m+2)}(x)  &\Delta w_1^{(m+2)}(x)   \\[2mm]
\end{bmatrix} 
\end{equation}

\medskip
\noindent Case II: If $n$ is even, i.e., $n=2m$ for some $m\in\mathbb{N}$, then we choose $m+2$ pairs of initial conditions $\{w_0^{(i)}, w_1^{(i)}\}$, $i=1, \dots, m+2$, and a boundary condition $h$ so that they satisfy (\ref{regularity}) and 
\begin{equation}\label{positivity2'}
|\det{\widetilde{W}(x)}|\geq r_0, \ a.e.\ x\in\Omega
\end{equation}
where $\widetilde{W}(x)$ is the $(n+3)\times(n+3)$ matrix defined by
\begin{equation}\label{10}
\widetilde{W}(x) = 
\begin{bmatrix}
w_0^{(1)}(x) & w_1^{(1)}(x) & \dis\pa_{x_1} w_0^{(1)}(x) & \cdots & \dis\pa_{x_n} w_0^{(1)}(x) &\Delta w_0^{(1)}(x) \\[2mm]
w_1^{(1)}(x) & w_{tt}^{(1)}(x) &\pa_{x_1}w_1^{(1)}(x)& \cdots &\pa_{x_n}w_1^{(1)}(x) &\Delta w_1^{(1)}(x) \\[2mm]
\vdots & \vdots & \vdots & \ddots & \vdots & \vdots \\
w_0^{(m+1)}(x) & w_1^{(m+1)}(x) &\dis\pa_{x_1} w_0^{(m+1)}(x) & \cdots & \dis\pa_{x_n} w_0^{(m+1)}(x) & \Delta w_0^{(m+1)}(x)\\[2mm]
 w_1^{(m+1)}(x) & w_{tt}^{(m+1)}(x) &\pa_{x_1}w_1^{(m+1)}(x)& \cdots &\pa_{x_n}w_1^{(m+1)}(x) &\Delta w_1^{(m+1)}(x) \\[2mm]
w_0^{(m+2)}(x) & w_1^{(m+2)}(x) &\dis\pa_{x_1} w_0^{(m+2)}(x) & \cdots & \dis\pa_{x_n} w_0^{(m+2)}(x) & \Delta w_0^{(m+2)}(x)
\end{bmatrix}
\end{equation}

Let $w^{(i)}(c,q_1,q_0,\q)$ and $w^{(i)}(\Tilde{c},p_1,p_0,\p)$ be the corresponding solutions of equation (\ref{1}) with different coefficients $\{c,q_1,q_0,\q\}$ and $\{\Tilde{c},p_1,p_0,\p\}$, as well as the initial and boundary conditions $\{w_0^{(i)}, w_1^{(i)}, h\}$, $i=1, \cdots, m+2$. If we have the same Neumann boundary traces over the observed part $\Gamma_1$ of the boundary and over the time interval $[-T, T]$, i.e., for $i = 1, \cdots, m+2$,
\begin{equation}\label{samemeasurement}
\frac{\pa w^{(i)}(c,q_1,q_0,\q)}{\pa\nu}(x,t) = \frac{\pa w^{(i)}(\Tilde{c},p_1,p_0,\p)}{\pa\nu}(x,t), \ (x,t) \in \Gamma_1\times [-T,T],
\end{equation}
then we must have that all the coefficients coincide, namely, 
\begin{equation}\label{13}
c (x) = \Tilde{c} (x), \ q_1 (x) = p_1 (x), \ q_0 (x) = p_0 (x), \ {\q}(x) = {\p}(x) \ \ a.e. \ x\in\Omega.
\end{equation}
\end{theorem}

After proving the above uniqueness theorem, we may also get the following Lipschitz stability result for recovering all coefficients $\{c,q_1,q_0,{\bf q}\}$ from the corresponding finite sets of boundary measurements.
\begin{theorem}\label{Th2}
Under the assumptions in Theorem~1.1, again let $w^{(i)}(c,q_1,q_0,\q)$ and $w^{(i)}(\Tilde{c},p_1,p_0,\p)$ denote the corresponding solutions of equation (\ref{1}) with coefficients $\{c,q_1,q_0,\q\}$ and $\{\Tilde{c},p_1,p_0,\p\}$, as well as the initial and boundary conditions $\{w_0^{(i)}, w_1^{(i)}, h\}$, $i=1, \cdots, m+2$ (either $n$ is odd or even). Then there exists $C>0$ depends on $\Omega$, $T$, $\Gamma_1$, $c$, $q_1$, $q_0$, $\q$, $w_0^{(i)}$, $w_1^{(i)}$, $h$ such that 
\begin{eqnarray}\label{stability}
& \ & \|c^2-{\Tilde{c}}^2\|_{L^2(\Omega)}^2 + \|q_1-p_1\|_{L^2(\Omega)}^2+\|q_0-p_0\|_{L^2(\Omega)}^2+\|{\q}-{\p}\|^2_{{\bf L}^2(\Omega)} \nonumber \\[2mm]
& \leq & C\sum_{i=1}^{m+2}\left\|\frac{\pa w_{tt}^{(i)}({c,q_1,q_0,\q})}{\pa\nu}-\frac{\pa w_{tt}^{(i)}({\Tilde{c},p_1,p_0,\p})}{\pa\nu}\right\|^2_{L^2(\Sigma_1)},
\end{eqnarray}
for all such coefficients $c, \tilde{c}, q_1,q_0, p_1, p_0\in H^1_0(\Omega)$, ${\q},{\p}\in \left(H^1_0(\Omega)\right)^n$, where $\|\cdot\|_{{\bf L}^2(\Omega)}$ is defined as 
\begin{equation*}
\|{\bf r}\|_{{\bf L}^2(\Omega)} = \displaystyle\left(\int_{\Omega}\sum_{i=1}^n|r_i(x)|^2\,dx\right)^{\frac{1}{2}}, \ \mbox{for} \ \ {\bf r}(x) = (r_1(x), \cdots, r_n(x)).
\end{equation*} 
\end{theorem}

\medskip
\noindent{\bf Inverse source problem}. The first step to solve the inverse problem above is to convert it into a corresponding inverse source problem. Indeed, if we let 
\be\label{relation}
\begin{split}
f_2(x)&=c^2(x)-\Tilde{c}^2(x),  \quad f_1(x)=p_1(x)-q_1(x),   \\[2mm]
f_0(x)&=p_0(x)-q_0(x),  \quad {\f}(x)={\p}(x)-{\q}(x); \\[2mm]
u(x,t)&=w({c,q_1,q_0,\q})-w({\Tilde{c},p_1,p_0,\p}), \ R(x,t)=w({\Tilde{c},p_1,p_0,\p})(x,t),
\end{split}
\ee
then $u=u(x,t)$ is readily seen to satisfy the following homogeneous mixed problem
\begin{equation}\label{2}
{
\begin{cases}
u_{tt} - c^2(x)\Delta u +q_1(x)u_t+q_0(x)u+{\q}(x)\cdot \nabla u 
=S(x,t)
& \mbox{in } Q\\[2mm]
u\left(x,0\right) = u_t\left(x,0\right) = 0 & \mbox{in } \Omega \\[2mm]
u(x,t) = 0 & \mbox{in } \Sigma,
\end{cases}
}
\end{equation}
where
\begin{equation}\label{rhs}
    S(x,t) = f_0(x)R(x,t)+f_1(x)R_t(x,t)+{\f}(x)\cdot \nabla R(x,t)+f_2(x)\Delta R(x,t).
\end{equation}
Here we assume that $c\in\mathcal{C}$, $q_0, q_1\in L^\infty(\Omega)$ and ${\q}\in \left(L^\infty(\Omega)\right)^n$ are given fixed and $R=R(x,t)$ is a given function that can be suitably chosen. On the other hand, the source coefficients $f_0,f_1,f_2\in L^2(\Omega)$ and ${\f}\in\left(L^2(\Omega)\right)^n$ are assumed to be unknown. The inverse source problem is to determine $f_0,f_1,f_2$ and ${\f}$ from the Neumann boundary measurements of $u$ over the observed part $\Gamma_1$ of the boundary and over a sufficiently long time interval $[-T, T]$.  More specifically, corresponding with Theorems \ref{Th1} and \ref{Th2}, we will prove the following uniqueness and stability results. 

\begin{theorem}\label{Th3}
Under geometrical assumptions (A.1) and (A.2) and let $T$ satisfy (\ref{largetime}). 
Depending on the dimension $n$, we assume the following regularity and positivity conditions:

\smallskip
Case I: If $n$ is odd, i.e., $n=2m+1$ for some $m\in\mathbb{N}$, then we choose $m+2$ functions $R^{(1)}$, $\cdots$, $R^{(m+2)}$ such that they satisfy 
\begin{equation}\label{regularity2}
R^{(i)},R_{t}^{(i)},R_{tt}^{(i)},R_{ttt}^{(i)}\in W^{2,\infty}(Q), \ i=1, \cdots, m+2
\end{equation}
and there exists $r_0>0$ such that 
\begin{equation}\label{positivity1}
|\det{U}(x)|\geq r_0, \ a.e. \ x\in\Omega 
\end{equation}
where $U(x)$ is the $(n+3)\times(n+3)$ matrix defined by
\begin{equation}\label{19}
U(x) =
\begin{bmatrix}
R^{(1)}(x,0) & R_t^{(1)}(x,0) & \dis\pa_{x_1} R^{(1)}(x,0) & \cdots & \dis\pa_{x_n} R^{(1)}(x,0) &\Delta R^{(1)}(x,0) \\[2mm]
R_t^{(1)}(x,0) & R_{tt}^{(1)}(x,0) &\pa_{x_1}R_t^{(1)}(x,0)& \cdots & \pa_{x_n}R_t^{(1)}(x,0) &\Delta R_t^{(1)}(x,0) \\[2mm]
\vdots & \vdots & \vdots & \ddots & \vdots & \vdots \\
R^{(m+2)}(x,0) & R_t^{(m+2)}(x,0) &\dis\pa_{x_1} R^{(m+2)}(x,0) & \cdots & \dis\pa_{x_n} R^{(m+2)}(x,0) & \Delta R^{(m+2)}(x,0)\\[2mm]
R_t^{(m+2)}(x,0) & R_{tt}^{(m+2)}(x,0) & \pa_{x_1}R_t^{(m+2)}(x,0)& \cdots & \pa_{x_n}R_t^{(m+2)}(x,0) &\Delta R_t^{(m+2)}(x,0) \\[2mm]
\end{bmatrix} 
\end{equation}

\medskip
Case II: If $n$ is even, i.e., $n=2m$ for some $m\in\mathbb{N}$, then we choose $m+2$ functions $R^{(1)}$, $\cdots$, $R^{(m+2)}$ such that they satisfy (\ref{regularity2}) and there exists $r_0>0$ such that
\begin{equation}\label{positivity2}
|\det{\widetilde{U}}(x)|\geq r_0, \ a.e. \ x\in\Omega 
\end{equation}
where $\widetilde{U}(x)$ is the $(n+3)\times(n+3)$ matrix defined by
\begin{equation}\label{22}
\widetilde{U}(x) =
\begin{bmatrix}
R^{(1)}(x,0) & R_t^{(1)}(x,0) & \dis\pa_{x_1} R^{(1)}(x,0) & \cdots & \dis\pa_{x_n} R^{(1)}(x,0) &\Delta R^{(1)}(x,0) \\[2mm]
R_t^{(1)}(x,0) & R_{tt}^{(1)}(x,0) &\pa_{x_1}R_t^{(1)}(x,0)& \cdots & \pa_{x_n}R_t^{(1)}(x,0) &\Delta R_t^{(1)}(x,0) \\[2mm]
\vdots & \vdots & \vdots & \ddots & \vdots & \vdots \\
R^{(m+1)}(x,0) & R_t^{(m+1)}(x,0) &\dis\pa_{x_1} R^{(m+1)}(x,0) & \cdots & \dis\pa_{x_n} R^{(m+1)}(x,0) & \Delta R^{(m+1)}(x,0)\\[2mm]
R_t^{(m+1)}(x,0) & R_{tt}^{(m+1)}(x,0) & \pa_{x_1}R_t^{(m+1)}(x,0)& \cdots & \pa_{x_n}R_t^{(m+1)}(x,0) &\Delta R_t^{(m+1)}(x,0) \\[2mm]
R^{(m+2)}(x,0) & R_t^{(m+2)}(x,0) &\dis\pa_{x_1} R^{(m+2)}(x,0) & \cdots & \dis\pa_{x_n} R^{(m+2)}(x,0) & \Delta R^{(m+2)}(x,0)\\[2mm]
\end{bmatrix} 
\end{equation}

Let $u^{(i)}(f_0,f_1,f_2,\f)$ be the solutions of equation (\ref{2}) with the functions $R^{(i)}$, $i=1, \cdots, m+2$. If
\begin{equation}\label{th5assume3}
\frac{\pa u^{(i)}({f_0,f_1,f_2,\f})}{\pa\nu}(x,t)=0, \ (x, t)\in \Gamma_1\times[-T,T], \ i=1, \cdots, m+2,
\end{equation}
then we must have
\begin{equation}\label{th5conclusion}
f_0(x)=f_1(x)=f_2(x)={\f}(x) = 0, \quad \mbox{ a.e. } x \in \Omega.
\end{equation}
\end{theorem}


\begin{theorem}\label{Th4}
Under the assumptions in Theorem \ref{Th3}, again let $u^{(i)}(f_0,f_1,f_2,\f)$ denote the solutions of equation (\ref{2}) with the functions $R^{(i)}$, $i=1, \cdots, m+2$ (either $n$ is odd or even). Then there exists $C>0$ depends on $\Omega$, $T$, $\Gamma_1$, $c$, $q_1$, $q_0$, $\q$, $w_0^{(i)}$, $w_1^{(i)}$, $h$ such that 
\begin{equation}\label{stability}
\|{f_0}\|^2_{L^2(\Omega)}+\|{f_1}\|^2_{L^2(\Omega)}+\|{f_2}\|^2_{L^2(\Omega)}+\|{\f}\|^2_{{\bf L}^2(\Omega)}\leq C\sum_{i=1}^{m+2}\left\|\frac{\pa u_{tt}^{(i)}({f_0,f_1,f_2,\f})}{\pa\nu}\right\|^2_{L^2(\Sigma_1)}
\end{equation}
for all $f_0, f_1, f_2\in H_0^1(\Omega)$ and ${\f} \in \left(H_0^{1}(\Omega)\right)^n$.
\end{theorem}

The rest of the paper is organized as follows. In the next section we recall some necessary tools to solve the inverse problem. This includes the sharp Carleman estimate, continuous observability inequality and regularity theory for general second-order hyperbolic equations with Dirichlet boundary condition. In Section 3 we provide the proofs of Theorems \ref{Th1}, \ref{Th2}, \ref{Th3} and \ref{Th4}, and in the last section we give some examples where the positivity conditions (\ref{positivity1'}), (\ref{positivity2'}), (\ref{positivity1}) and (\ref{positivity2}) are satisfied and some concluding remarks.

\section{Carleman Estimate, Continuous Observability Inequality and Regularity Theory for Second-Order Hyperbolic Equations}

In this section we recall some key ingredients of the proofs used in the next section. This includes Carleman estimate, continuous observability inequality, as well as regularity theory for general second-order hyperbolic equations with Dirichlet boundary condition. For simplicity here we only state the main results and refer to \cite{TY2002} and \cite{LLT1986} for greater details.

To begin with, consider a Riemannian metric $g(\cdot,\cdot)=\langle \cdot, \cdot \rangle$ and squared norm $|X|^2=g(X,X),$ on a smooth finite dimensional manifold $\mathcal{M}$. On the Riemannian manifold $(\mathcal{M}, g)$ we define $\Omega$ as an open bounded, connected set of $\mathcal{M}$ with smooth boundary $\Gamma=\overline{\Gamma_0\cup\Gamma_1}$, where $\Gamma_0\cap\Gamma_1=\emptyset$.  Let $\nu$ denote the unit outward normal field along the boundary $\Gamma$.  Furthermore, we denote by $\Delta_g$ the Laplace--Beltrami operator on the manifold $\mathcal{M}$ and by $D$ the Levi--Civita connection on $\mathcal{M}$.  

Consider the following second-order hyperbolic equation with energy level terms defined on $Q = \Omega\times[-T,T]$ for some $T>0$:
\be\label{mainequation}
w_{tt}(x,t) - \Delta_g w(x,t) + F(w) = G(x,t), \quad (x, t)\in Q=\Omega\times[-T,T]
\ee
where the forcing term $G\in L^2(Q)$ and the energy level differential term $F(w)$ is given by
\begin{equation*}\label{fw}
F(w)= \langle{\bf P}(x,t), Dw\rangle +P_1(x,t)w_t+P_0(x,t)w.
\end{equation*}
Here $P_0$, $P_1$ are functions on $\Omega\times[-T,T]$, ${\bf P}(x, t)$ is a vector field on $\mathcal{M}$ for $t\in[-T,T]$, and they satisfy the following estimate: there exists a constant $C_T>0$ such that 
\begin{equation*}\label{fwproperty}
|F(w)|\leq C_T[w^2+w_t^2+|Dw|^2], \ \forall(x, t)\in Q.
\end{equation*}

\noindent\textbf{Pseudo-convex function.} Having chosen, on the strength of geometrical assumption (A.1), a strictly convex function $d(x)$, we can define the function $\varphi(x,t): \Omega\times\mathbb{R}\to\mathbb{R}$ of class $C^3$
by setting
\begin{equation*}\label{pseudo1}
\dis\varphi(x,t)=d(x)-\alpha t^2, \quad x\in\Omega, \  t\in[-T,T],
\end{equation*}
where $T>T_0$ as in (\ref{largetime}). Moreover, $\alpha\in(0,1)$ is selected as follows: Let $T>T_0$ be given, then there exists $\delta>0$ such that
\begin{equation*}\label{pseudo3}
\dis T^2>4\max_{x\in\overline{\Omega}}d(x)+4\delta.
\end{equation*}
For this $\delta>0$, there exists a constant $\alpha\in(0,1)$, such that
\begin{equation*}\label{pseudo4}
\dis \alpha T^2>4\max_{x\in\overline{\Omega}}d(x)+4\delta.
\end{equation*}
It is easy to check such function $\varphi(x,t)$ carries the following properties:

\medskip
(a) For the constant $\delta>0$ fixed above, we have 
\begin{equation*}\label{pseudoproperty1}
\dis\varphi(x,-T)=\varphi(x,T) \leq\max_{x\in\overline{\Omega}}d(x)-\alpha T^2\leq-\delta
\ \text{uniformly in} \ x\in\Omega;
\end{equation*}
and
\begin{equation*}\label{pseudoproperty2}
\dis\varphi(x,t)\leq\varphi\left(x,0\right), \quad \text{for any} \ t\in[-T,T] \ \text{and any} \ x\in\Omega.
\end{equation*}

(b) There are $t_0$ and $t_1$, with $-T<t_0<0<t_1<T$, say, chosen symmetrically about $0$, such that
\begin{equation*}\label{assume4}
\min_{x\in\overline{\Omega},t\in[t_0,t_1]}\varphi(x,t)\geq\sigma,
\quad \text{where} \ 0<\sigma<m_0 = \min_{x\in\overline{\Omega}}d(x).
\end{equation*}
Moreover, let $Q(\sigma)$ be the subset of $Q=\Omega\times[-T,T]$ defined by
\begin{equation}\label{qsigma}
Q{(\sigma)}=\{(x,t): \varphi(x,t)\geq\sigma>0, x\in\Omega, -T\leq t\leq T\},
\end{equation}
Then we have
\be\label{qsigmaproperty}
\Omega\times[t_0,t_1]\subset Q(\sigma)\subset\Omega\times[-T,T].
\ee

\noindent\textbf{Carleman estimate for general second-order hyperbolic equations}. 
We now return to the equation \eqref{mainequation}, and consider solutions $w(x,t)$ in the class
\begin{equation}\label{h1}
\begin{cases}
w\in H^{1,1}(Q) = L^2(-T,T;H^1(\Omega))\cap H^1(-T,T;L^2(\Omega)); \\[2mm]
w_t, \frac{\pa w}{\pa\nu}\in L^2(-T, T; L^2(\Gamma)).
\end{cases}
\end{equation}

Then for these solutions with geometrical assumptions (A.1) and (A.2) on $\Omega$, the following one-parameter family of estimates hold true, with $\beta>0$ being a suitable constant ($\beta$ is positive by virtue of (A.2)), for all $\tau>0$ sufficiently large and $\epsilon>0$ small:
\begin{multline}\label{carleman}
BT(w)+2\int_{Q}e^{2\tau\varphi}|G|^2\,dQ+C_{1,T}e^{2\tau\sigma}\int_{Q}w^2\,dQ+ c_T\tau^3e^{-2\tau\delta}[E_w(-T)+E_w(T)] \\[2mm]
\geq C_{1,\tau}\int_Qe^{2\tau\varphi}[w_t^2+|Dw|^2]\,dQ + C_{2,\tau}\int_{Q{(\sigma)}}e^{2\tau\varphi}w^2\,dxdt
\end{multline}
where
\begin{equation}\label{c1tauc2tau}
C_{1,\tau}=\tau\epsilon(1-\alpha)-2C_T, \quad C_{2,\tau}=2\tau^3\beta+\mathcal{O}(\tau^2)-2C_T.
\end{equation}
Here $\delta>0$, $\sigma>0$ are the constants as in above, $C_T$, $c_T$ and $C_{1,T}$ are positive constants depending on $T$, as well as $d$ (but not on $\tau$).
The energy function $E_w(t)$ is defined as
\begin{equation*}\label{energy}
E_w(t)=\int_{\Omega}[w^2(x,t)+w_t^2(x,t)+|Dw(x,t)|^2]\,d\Omega.
\end{equation*}
In addition, $BT(w)$ stands for boundary terms and can be explicitly calculated as
\begin{eqnarray*}\label{boundary}
BT(w) & = & 2\tau\int_{\Sigma}e^{2\tau\varphi}\left(w_t^2-|Dw|^2\right)\langle Dd,\nu\rangle\,d\Sigma \nonumber\\[2mm]
&+& 4\tau\int_{\Sigma}e^{2\tau\varphi}\langle Dd, Dw\rangle\langle Dw,\nu\rangle\,d\Sigma + 8\alpha\tau\int_{\Sigma}e^{2\tau\varphi}tw_t \langle Dw,\nu\rangle\,d\Sigma \nonumber\\[2mm]
&+& 4\tau^2\int_{\Sigma}e^{2\tau\varphi}\left[|Dd|^2-4\alpha^2t^2+\frac{\Delta d-\alpha-1}{2\tau}\right]w \langle Dw,\nu\rangle\,d\Sigma \nonumber \\[2mm]
&+& 2\tau\int_{\Sigma}e^{2\tau\varphi}\left[2\tau^2\left(|Dd|^2-4\alpha^2t^2\right)+ \tau(3\alpha+1)\right]w^2\langle Dd,\nu\rangle\,d\Sigma.
\end{eqnarray*}
Clearly if we have $w|_{\Gamma\times[-T,T]} = 0$ and  $\frac{\pa w}{\pa\nu}=\langle Dw, \nu\rangle = 0$ on $\Gamma_1\times[-T,T]$,  
then in view of the geometrical assumption (A.1) we may compute 
\begin{equation}\label{boundary}
BT(w) = 2\tau\int_{-T}^{T}\int_{\Gamma_0}e^{2\tau\varphi}|D w|^2 \langle Dd,\nu\rangle\,d\Gamma_0{d}t \leq 0. 
\end{equation} 

\medskip
\noindent\textbf{Continuous observability inequality}. As a corollary of the Carleman estimate, we also have the following continuous observability inequality 
\be\label{coi}
C_T E_w(0)\leq\int_{-T}^T\int_{\Gamma_1} \left(\frac{\pa w}{\pa\nu}\right)^2 d\Gamma dt+\|G\|^2_{L^2(Q)}
\ee
for the equation (\ref{mainequation}) with homogeneous Dirichlet boundary condition $w|_{\Sigma}=0$. Here $T>T_0$ as in \eqref{largetime} and $\Omega$ satisfies the geometrical assumptions (A.1) and (A.2).
\begin{remark}
The continuous observability inequality (\ref{coi}) may be interpreted as follows: If the second-order hyperbolic equation equation (\ref{mainequation}) has homogeneous Dirichlet boundary condition and nonhomogeneous forcing term $G\in L^2(Q)$, and Neumann boundary trace $\frac{\pa w}{\pa\nu}\in L^2(\Sigma_1)$, then necessarily the initial conditions $\{w(\cd, 0), w_t(\cd,0)\}$ must lie in the natural energy space $H^1_0(\Omega)\times L^2(\Omega)$. This fact will be used in the proofs in Section 3.
\end{remark}

\noindent\textbf{Regularity theory for general second-order hyperbolic equations with Dirichlet boundary condition}. Consider the second-order hyperbolic equation (\ref{mainequation}) with initial conditions $w(x, 0)=w_0(x)$, $w_t(x,0)=w_1(x)$ and Dirichlet boundary condition $w|_{\Sigma} = h(x,t)$. Then the following interior and boundary regularity results for the solution $w$ hold true:
For $\gamma\geq 0$ (not necessarily an integer), if the given data satisfy the following regularity assumptions 
\begin{equation*}
\begin{cases}
G\in L^1(0,T; H^{\gamma}(\Omega)), \ \pa_t^{(\gamma)}G\in L^1(0,T; L^2(\Omega)), \\[2mm]
w_0\in H^{{\gamma}+1}(\Omega), \ w_1\in H^{\gamma}(\Omega), \ h\in H^{{\gamma}+1}(\Sigma)
\end{cases}
\end{equation*}
with all compatibility conditions (trace coincidence) which make sense. Then, we have the following regularity for the solution $w$:
\be\label{reg}
w\in C([0,T]; H^{{\gamma}+1}(\Omega)), \ \pa_t^{({\gamma}+1)}w\in C([0,T]; L^2(\Omega)); \ \frac{\pa w}{\pa\nu}\in H^{\gamma}(\Sigma).
\ee

\section{Main Proofs}

In this section we give the main proofs of the uniqueness and stability results established in the first section. We focus on proving Theorems \ref{Th3} and \ref{Th4} for the inverse source problem since Theorems \ref{Th1} and \ref{Th2} of the original inverse problem will then follow from the relation \eqref{relation} between the two problems and the regularity theory result recalled in Section 2. 
Henceforth for convenience we use $C$ to denote a generic positive constant which may depend on $\Omega$, $T$, $c$, $q_1$, $q_0$, ${\q}$, $r_0$, $w^{(i)}$, $u^{(i)}$, $R^{(i)}$, $i=1, \cdots, m+2$, but not on the free large parameter $\tau$ appearing in the Carleman estimate.

\medskip 
\noindent\textbf{Proof of Theorem~\ref{Th3}}. 
First we consider the case when $n$ is odd, i.e., $n=2m+1$, for some $m\in\mathbb{N}$. Then corresponding with the choice of $R^{(i)}$, $i=1, \cdots, m+2$, we have $m+2$ equations of the form \eqref{2} with solutions $u^{(i)}=u^{(i)}(x,t)$ that satisfy 
\begin{equation}\label{th1u}
\begin{cases}
u^{(i)}_{tt} - c^2(x)\Delta u^{(i)} + q_1(x)u_t^{(i)}+q_0(x)u^{(i)}+{\q}(x)\cdot\nabla u^{(i)} = S^{(i)}(x,t) & \mbox{in } \ Q \\[2mm]
u^{(i)}(x, 0)  = u^{(i)}_t(x, 0) = 0 & \mbox{in } \ \Omega  \\[2mm]
u^{(i)}|_{\Gamma\times [-T,T]}=0, \quad \frac{\pa u^{(i)}}{\pa\nu}|_{\Gamma_1\times[-T,T]}=0 & \mbox{in } \Sigma, \Sigma_1,
\end{cases}
\end{equation}
where $S^{(i)}(x,t)$ is defined in (\ref{rhs}) with $R$ being replaced by $R^{(i)}$. 

Note since $c\in\mathcal{C}$, $q_1, q_0\in L^{\infty}(\Omega)$ and ${\bf q}\in\left(L^{\infty}(\Omega)\right)^n$,
the equation in (\ref{th1u}) can be written as a Riemannian wave equation with respect to the metric $g={c^{-2}(x)}dx^2$, modulo lower-order terms\footnote{More precisely, we have $\disp \Delta_gu = c^2\Delta u + c^n\nabla(c^{2-n})\cdot\nabla u$} 
\begin{equation*}
u^{(i)}_{tt} - \Delta_g u^{(i)} + \mbox{``lower-order terms''} = S^{(i)}(x,t).
\end{equation*}
Moreover, by the regularity assumption \eqref{regularity2}, we have that $S^{(i)}\in L^2(Q)$ and by Cauchy--Schwarz inequality
\begin{equation*}
|S^{(i)}(x,t)|^2 \leq C\left(|f_0(x)|^2+|f_1(x)|^2+|{\f}(x)|^2+|f_2(x)|^2\right).
\end{equation*}
Thus we  can apply the Carleman estimate \eqref{carleman} for solution $u^{(i)}$ in the class \eqref{h1} and get the following inequality for sufficiently large $\tau$:
\begin{equation}\label{th1carleman1}
\begin{split}
\tau\int_{Q}e^{2\tau\varphi}[(u^{(i)}_t)^2&+|D u^{(i)}|^2]dQ+\tau^3\int_{Q{(\sigma)}}e^{2\tau\varphi}(u^{(i)})^2dxdt\\
&\leq C\int_{Q}e^{2\tau\varphi}\left(|f_0(x)|^2+|f_1(x)|^2+|{\f}(x)|^2+|f_2(x)|^2\right) dQ+Ce^{2\tau\sigma}.
\end{split}
\end{equation}
Note here we have dropped the unnecessary terms in the Carleman estimate (\ref{carleman}) as well as the boundary terms $BT(u^{(i)})$ since the homogeneous boundary data $u^{(i)}|_{\Gamma\times[-T,T]}=\frac{\pa u^{(i)}}{\pa\nu}|_{\Gamma_1\times[-T,T]}=0$ imply $BT(u^{(i)})\leq 0$, as suggested in (\ref{boundary}).

Differentiate the $u^{(i)}$-system \eqref{th1u} in time $t$, we get the following $u^{(i)}_t$-problem
\begin{equation}\label{th1ut}
\begin{cases}
(u^{(i)}_t)_{tt} - c^2(x)\Delta u^{(i)}_t +q_1(x)(u_t^{(i)})_t+q_0(x)u^{(i)}_t+{\q}(x)\cdot\nabla u^{(i)}_t=S^{(i)}_t(x,t) & \mbox{in } \ Q \\[2mm]
(u_t^{(i)})(x, 0)=0, \ (u^{(i)}_t)_t(x, 0) = S^{(i)}(x,0) & \mbox{in } \ \Omega  \\[2mm]
u^{(i)}_t|_{\Gamma\times [-T,T]}=0, \quad \frac{\pa u^{(i)}_t}{\pa\nu}|_{\Gamma_1\times[-T,T]}=0 & \mbox{in } \Sigma, \Sigma_1.
\end{cases}
\end{equation}
Note again by \eqref{regularity2} we have $S^{(i)}_t\in L^2(Q)$ and by Cauchy--Schwarz inequality
\begin{equation*}
|S^{(i)}_t(x,t)|^2 \leq C\left(|f_0(x)|^2+|f_1(x)|^2+|{\f}(x)|^2+|f_2(x)|^2\right).
\end{equation*}
In addition, $BT(u_t^{(i)})\leq 0$ since $u_t^{(i)}|_{\Gamma\times[-T,T]}=\frac{\pa u_t^{(i)}}{\pa\nu}|_{\Gamma_1\times[-T,T]}=0$. 
Thus similar to (\ref{th1carleman1}) we can apply Carleman estimate \eqref{carleman} for solutions $u_t^{(i)}$ in the class \eqref{h1} and get the following inequality for sufficiently large $\tau$:
\be\label{th1carleman2}
\begin{split}
\tau\int_{Q}e^{2\tau\varphi}[(u^{(i)}_{tt})^2&+|D u^{(i)}_t|^2]dQ+\tau^3\int_{Q{(\sigma)}}e^{2\tau\varphi}(u^{(i)}_t)^2dxdt\\
&\leq C\int_{Q}e^{2\tau\varphi}\left(|f_0(x)|^2+|f_1(x)|^2+|{\f}(x)|^2+|f_2(x)|^2\right) dQ+Ce^{2\tau\sigma}.
\end{split}
\ee

Continue with this process, we differentiate \eqref{th1ut} in $t$ two more times, and get the corresponding $u^{(i)}_{tt}$ and $u^{(i)}_{ttt}$-systems
\begin{equation}\label{th1utt}
\begin{cases}
(u^{(i)}_{tt})_{tt} - c^2(x) \Delta u^{(i)}_{tt}+q_1(x)(u^{(i)}_{tt})_t+q_0(x)u^{(i)}_{tt}+ {\q}(x)\cdot\nabla u^{(i)}_{tt} = S^{(i)}_{tt}(x,t)\\[2mm]
u^{(i)}_{tt}(x, 0)  = S^{(i)}(x,0), \ (u^{(i)}_{tt})_t(x, 0) = S^{(i)}_t(x,0)-q_1(x)S^{(i)}(x,0)  \\[2mm]
u^{(i)}_{tt}|_{\Gamma\times [-T, T]}=0, \quad \frac{\pa u^{(i)}_{tt}}{\pa\nu}|_{\Gamma_1\times[-T,T]}=0
\end{cases}
\end{equation}
and
\begin{equation}\label{th1uttt}
\begin{cases}
(u^{(i)}_{ttt})_{tt} - c^2(x)\Delta u^{(i)}_{ttt}+ q_1(x)(u^{(i)}_{ttt})_t+ q_0(x)u^{(i)}_{ttt}+{\q}(x)\cdot\nabla u^{(i)}_{ttt} = S^{(i)}_{ttt}(x,t)\\[2mm]
(u^{(i)}_{ttt})(x, 0) = S^{(i)}_t(x,0) -q_1(x)S^{(i)}(x,0)\\[2mm] 
(u^{(i)}_{ttt})_t(x, 0) = S^{(i)}_{tt}(x,0) + c^2\Delta S^{(i)}(x,0)-q_1S^{(i)}_t(x,0) - q_0S^{(i)}(x,0)-{\bf q}\cdot \nabla S^{(i)}(x,0) \\[2mm]
u^{(i)}_{ttt}|_{\Gamma\times [-T, T]}=0, \quad \frac{\pa u^{(i)}_{ttt}}{\pa\nu}|_{\Gamma_1\times[-T,T]}=0.
\end{cases}
\end{equation}

Again by \eqref{regularity2}, Cauchy--Schwarz inequality and the homogeneous Dirichlet and Neumann boundary data, we can apply Carleman estimate \eqref{carleman} to the corresponding $u^{(i)}_{tt}$, $u^{(i)}_{ttt}$-systems above and get the following inequalities that are similar to \eqref{th1carleman1} and \eqref{th1carleman2}, for $\tau$ sufficiently large

\be\label{th1carleman3}
\begin{split}
\tau\int_{Q}e^{2\tau\varphi}[(u^{(i)}_{ttt})^2&+|Du^{(i)}_{tt}|^2]dQ+\tau^3\int_{Q{(\sigma)}}e^{2\tau\varphi}(u^{(i)}_{tt})^2dxdt\\
&\leq C\int_{Q}e^{2\tau\varphi}\parens{|f_0(x)|^2+|f_1(x)|^2+|{\f}(x)|^2+|f_2(x)|^2} dQ+Ce^{2\tau\sigma}.
\end{split}
\ee

\be\label{th1carleman3'}
\begin{split}
\tau\int_{Q}e^{2\tau\varphi}[(u^{(i)}_{tttt})^2&+|D u^{(i)}_{ttt}|^2]dQ+\tau^3\int_{Q{(\sigma)}}e^{2\tau\varphi}(u^{(i)}_{ttt})^2dxdt\\
&\leq C\int_{Q}e^{2\tau\varphi}\parens{|f_0(x)|^2+|f_1(x)|^2+|{\f}(x)|^2+|f_2(x)|^2} dQ+Ce^{2\tau\sigma}.
\end{split}
\ee

Add the four inequalities \eqref{th1carleman1}, \eqref{th1carleman2}, \eqref{th1carleman3}, \eqref{th1carleman3'} together, we get
\begin{equation}\label{th1carleman4}
\begin{split}
\tau&\int_{Q}e^{2\tau\varphi}[(u^{(i)}_{tttt})^2+(u^{(i)}_{ttt})^2+(u^{(i)}_{tt})^2+(u^{(i)}_t)^2+|D u^{(i)}_{ttt}|^2+|D u^{(i)}_{tt}|^2+|D u^{(i)}_{t}|^2+|D u^{(i)}|^2]dQ \\
&+ \tau^3\int_{Q{(\sigma)}}e^{2\tau\varphi}[(u^{(i)}_{ttt})^2+(u^{(i)}_{tt})^2+(u^{(i)}_t)^2+(u^{(i)})^2]dxdt\\
&\leq C\int_{Q}e^{2\tau\varphi}\parens{|f_0(x)|^2+|f_1(x)|^2+|{\f}(x)|^2+|f_2(x)|^2} dQ+Ce^{2\tau\sigma}.
 \end{split}
\end{equation}

We now analyze the integral term on the right-hand side of \eqref{th1carleman4}. First note that by estimating the $u^{(i)}$-equation in (\ref{th1u}) and $u^{(i)}_t$-equation in (\ref{th1ut}) at time $t=0$, we can get
\begin{equation}\label{th1eq1}
\begin{cases}
u^{(i)}_{tt}(x,0) = S^{(i)}(x,0) \\[2mm]
u^{(i)}_{ttt}(x,0) = S^{(i)}_t(x,0)-q_1(x)S^{(i)}(x,0).
\end{cases}
\end{equation}

Note the above equations hold for any $i$, $1\leq i\leq m+2$, so putting all of them together we get a $(n+3)\times (n+3)$ linear system
\begin{equation}\label{linearsystem}
\left[u^{(1)}_{tt}(x,0), u^{(1)}_{ttt}(x,0),\cdots,u^{(m+2)}_{tt}(x,0), u^{(m+2)}_{ttt}(x,0)\right]^T =
U_{q_1}(x)\left[f_0(x), f_1(x), {\f}(x), f_2(x)\right]^T
\end{equation}
where the coefficient matrix $U_{q_1}(x)$ is defined as 
\begin{equation}\label{uq}
U_{q_1}(x) =
\begin{bmatrix}
R^{(1)}(x,0) & R_t^{(1)}(x,0) & \dis\pa_{x_1} R^{(1)}(x,0) & \cdots & \dis\pa_{x_n} R^{(1)}(x,0) &\Delta R^{(1)}(x,0) \\[2mm]
\tilde{a}^{(1)}(x) & \tilde{b}^{(1)}(x) &\tilde{m}_1^{(1)}(x)& \cdots & \tilde{m}_n^{(1)}(x) &\tilde{\ell}^{(1)}(x) \\[2mm]
\vdots & \vdots & \vdots & \ddots & \vdots & \vdots \\
R^{(m+2)}(x,0) & R_t^{(m+2)}(x,0) &\dis\pa_{x_1} R^{(m+2)}(x,0) & \cdots & \dis\pa_{x_n} R^{(m+2)}(x,0) & \Delta R^{(m+2)}(x,0)\\[2mm]
\tilde{a}^{(m+2)}(x) & \tilde{b}^{(m+2)}(x) & \tilde{m}_1^{(m+2)}(x)& \cdots & \tilde{m}_n^{(m+2)}(x) &\tilde{\ell}^{(m+2)}(x) \\[2mm]
\end{bmatrix} 
\end{equation}
with
\begin{multline}\label{rndcoef}
    \tilde{a}^{(i)}(x) =R_t^{(i)}(x,0)-q_1(x)R^{(i)}(x,0), \ \ \tilde{b}^{(i)}(x) =R_{tt}^{(i)}(x,0)-q_1(x)R_t^{(i)}(x,0), \\[2mm]
    \tilde{m}_k^{(i)}(x) =\dis\pa_{x_k}R_t^{(i)}(x,0)-q_1(x)\dis\pa_{x_k}R^{(i)}(x), \ \ \tilde{\ell}^{(i)}(x,0) =\Delta R_t^{(i)}(x,0)-q_1(x)\Delta R^{(i)}(x,0).
\end{multline}
Notice that from doing elementary row operations, specifically, adding $q_1$ multiplied by an odd row to the subsequent even row, the matrix $U_{q_1}(x)$ and $U(x)$ as defined in (\ref{19}) have the same determinant. Thus the positivity assumption (\ref{positivity1}) implies that we may invert $U_{q_1}(x)$ in (\ref{uq}) to obtain
\begin{equation}\label{key}
\begin{split}
|f_0(x)|^2+|f_1(x)|^2+|f_2(x)|^2+|{\f}(x)|^2 & \leq C \sum_{i=1}^{m+2}\left(|u^{(i)}_{tt}(x,0)|^2+|u^{(i)}_{ttt}(x,0)|^2\right)\\[2mm]
& = C\left(|{\bf u}_{tt}(x,0)|^2+|{\bf u}_{ttt}(x,0)|^2\right)
\end{split}
\end{equation}
where we denote ${\bf u}(x,t) = (u^{(1)}(x,t), u^{(2)}(x,t), \cdots, u^{(m+2)}(x,t))$.
Thus by properties of the pseudo-convex function $\varphi$ and the Cauchy--Schwarz inequality, we can get the following estimate
\begin{eqnarray}\label{th1mainineq}
& \ & \int_{Q}e^{2\tau\varphi(x,t)}\parens{|f_0(x)|^2+|f_1(x)|^2+|{\f}(x)|^2+|f_2(x)|^2}dQ \\
               & \leq & C\int_{\Omega}\int_{-T}^Te^{2\tau\varphi(x,0)}\left(|{\bf u}_{tt}(x,0)|^2+|{\bf u}_{ttt}(x,0)|^2\right)dt\,d\Omega  \nonumber\\
                             & \leq &  C\left(\int_{\Omega}\int_{-T}^{0} \frac{d}{ds}[e^{2\tau\varphi(x,s)}\left(|{\bf u}_{tt}(x,s)|^2+|{\bf u}_{ttt}(x,s)|^2\right)]ds\,d\Omega \right. \nonumber\\
                             & \ & \ + \left. \int_{\Omega}e^{2\tau\varphi(x,-T)}\left(|{\bf u}_{tt}(x,-T)|^2+|{\bf u}_{ttt}(x,-T)|^2\right)d\Omega \right) \nonumber\\
                             & \leq &  C\left(\tau\int_{\Omega}\int_{-T}^{0}e^{2\tau\varphi(x,s)}\left(|{\bf u}_{tt}(x,s)|^2+|{\bf u}_{ttt}(x,s)|^2\right)]ds\,d\Omega\right. \nonumber\\
                             & \ & \  + \left.2\int_{\Omega}\int_{-T}^{0}e^{2\tau\varphi(x,s)}\left({\bf u}_{tt}\cdot {\bf u}_{ttt}+{\bf u}_{ttt}\cdot {\bf u}_{tttt}\right)]ds\,d\Omega \right. \nonumber\\
                             & \ & \ + \left. \int_{\Omega}e^{2\tau\varphi(x,-T)}\left(|{\bf u}_{tt}(x,-T)|^2+|{\bf u}_{ttt}(x,-T)|^2\right)d\Omega\right) \nonumber\\
                             & \leq & C\left(\tau\int_{Q} e^{2\tau\varphi}|{\bf u}_{tt}|^2dQ+\tau\int_{Q}e^{2\tau\varphi}|{\bf u}_{ttt}|^2dQ+\int_{Q}e^{2\tau\varphi}|{\bf u}_{tttt}|^2dQ\right).\nonumber
\end{eqnarray} 

Taking \eqref{th1mainineq} into \eqref{th1carleman4}, and note \eqref{th1carleman4} holds for all $i=1, \cdots m+2$, thus summing over $i$ in \eqref{th1carleman4} and dropping the non-negative gradient terms on the left-hand side, we get that for $\tau$ sufficiently large
\begin{eqnarray}\label{th1eq3}
& \ &\tau\int_{Q}e^{2\tau\varphi}\left(|{\bf u}_{tttt}|^2+|{\bf u}_{ttt}|^2+|{\bf u}_{tt}|^2+|{\bf u}_{t}|^2\right)\,dQ \\
& + & \tau^3\int_{Q{(\sigma)}}e^{2\tau\varphi}\left(|{\bf u}_{ttt}|^2+|{\bf u}_{tt}|^2+|{\bf u}_{t}|^2+|{\bf u}|^2\right)\,dx\,dt\nonumber\\[2mm]
& \leq & C\tau\int_{Q} e^{2\tau\varphi}(|{\bf u}_{tt}|^2+|{\bf u}_{ttt}|^2)dQ+C\int_{Q}e^{2\tau\varphi}|{\bf u}_{tttt}|^2dQ+Ce^{2\tau\sigma}. \nonumber
\end{eqnarray}

Since $e^{2\tau\varphi} < e^{2\tau\sigma}$ on $Q\backslash Q(\sigma)$ from the definition of $Q(\sigma)$ (\ref{qsigma}), we have the following
\begin{eqnarray*}\label{th1eq4}
& \ & \int_Q e^{2\tau\varphi}\left(|{\bf u}_{tt}|^2+|{\bf u}_{ttt}|^2\right)dQ \\
& = & \int_{Q(\sigma)} e^{2\tau\varphi}\left(|{\bf u}_{tt}|^2+|{\bf u}_{ttt}|^2\right)dx \, dt +\int_{Q\backslash Q(\sigma)} e^{2\tau\varphi}\left(|{\bf u}_{tt}|^2+|{\bf u}_{ttt}|^2\right)dx\, dt
\nonumber \\
& \leq & \int_{Q(\sigma)} e^{2\tau\varphi}\left(|{\bf u}_{tt}|^2+|{\bf u}_{ttt}|^2\right)dx \, dt
+ \  e^{2\tau\sigma} \int_{Q\backslash Q(\sigma)}\left(|{\bf u}_{tt}|^2+|{\bf u}_{ttt}|^2\right) dx \, dt
\end{eqnarray*}
and therefore \eqref{th1eq3} becomes
\begin{eqnarray}\label{th1eq5}
& \ &\tau\int_{Q}e^{2\tau\varphi}\left(|{\bf u}_{tttt}|^2+|{\bf u}_{ttt}|^2+|{\bf u}_{tt}|^2+|{\bf u}_{t}|^2\right)dQ \\
& +& \tau^3\int_{Q{(\sigma)}}e^{2\tau\varphi}\left(|{\bf u}_{ttt}|^2+|{\bf u}_{tt}|^2+|{\bf u}_{t}|^2+|{\bf u}|^2\right)dx\, dt\nonumber\\[2mm]
&\leq& C\tau\int_{Q(\sigma)} e^{2\tau\varphi}\left(|{\bf u}_{tt}|^2+|{\bf u}_{ttt}|^2\right)dx \, dt+C\int_{Q}e^{2\tau\varphi}|{\bf u}_{tttt}|^2dQ+Ce^{2\tau\sigma}.\nonumber
\end{eqnarray}

Note that in \eqref{th1eq5} the first and second terms on the right-hand side can be absorbed by the corresponding terms on the left-hand side when $\tau$ is taken large enough. Hence we may get the following estimate for sufficiently large $\tau$: 
\begin{equation*}\label{th1eq6}
\tau^3\int_{Q{(\sigma)}}e^{2\tau\varphi}\left(|{\bf u}_{ttt}|^2+|{\bf u}_{tt}|^2+|{\bf u}_{t}|^2+|{\bf u}|^2\right)dx \, dt\leq C\tau e^{2\tau\sigma}.
\end{equation*}
Use again the fact that $\varphi(x,t)\geq\sigma$ on $Q(\sigma)$ we hence get
\begin{equation*}\label{th1eq7}
\tau^2\int_{Q{(\sigma)}}|{\bf u}_{ttt}|^2+|{\bf u}_{tt}|^2+|{\bf u}_{t}|^2+|{\bf u}|^2 \ dx \, dt\leq C.
\end{equation*}

Since $\tau>0$ in a free large parameter and the constants $C$ do not depend on $\tau$, the above inequality implies we must have ${\bf u}={\bf 0}$ a.e. on $Q(\sigma)$.
Note from (\ref{qsigmaproperty}) the subspace $Q(\sigma)$ satisfies the property $\Omega\times[t_0, t_1]\subset Q(\sigma)\subset Q$ with $t_0<0<t_1$, therefore by evaluating the ${\bf u}$ and ${\bf u}_t$-systems of equations at $t=0$, we get the $(n+3)\times (n+3)$ linear system (see (\ref{linearsystem}))
\begin{equation*}
U_{q_1}(x) [f_0(x), f_1(x), {\bf f}(x), f_2(x)]^T = {\bf 0}, \ a.e. \ x\in\Omega.
\end{equation*}
As the coefficient matrix $U_{q_1}(x)$ is invertible from assumption (\ref{positivity1}), we must have the desired conclusion
\begin{equation*}
f_0(x)=f_1(x)=f_2(x)={\bf f}(x)= 0, \ a.e. \ x\in\Omega. 
\end{equation*}

For the case when $n$ is even, i.e., $n=2m$, $m\in\mathbb{N}$, we can basically repeat the above proof with obvious adjustments. The only difference here is that since $n=2m$ is even, the linear system (\ref{linearsystem}) contains an odd number ($n+3$) of equations. Therefore we only need $m+1$ pairs of equations from (\ref{th1eq1}) plus one more equation from $u_{tt}^{(m+2)}(x,0)$. Doing  this yields the matrix $\widetilde{U}_{q_1}(x)$, where
\begin{equation}\label{22}
\widetilde{U}_{q_1}(x) =
\begin{bmatrix}
R^{(1)}(x,0) & R_t^{(1)}(x,0) & \dis\pa_{x_1} R^{(1)}(x,0) & \cdots & \dis\pa_{x_n} R^{(1)}(x,0) &\Delta R^{(1)}(x,0) \\[2mm]
\tilde{a}^{(1)}(x) & \tilde{b}^{(1)}(x) &\tilde{m}_1^{(1)}(x)& \cdots & \tilde{m}_n^{(1)}(x) &\tilde{\ell}^{(1)}(x) \\[2mm]
\vdots & \vdots & \vdots & \ddots & \vdots & \vdots \\
R^{(m+1)}(x,0) & R_t^{(m+1)}(x,0) &\dis\pa_{x_1} R^{(m+1)}(x,0) & \cdots & \dis\pa_{x_n} R^{(m+1)}(x,0) & \Delta R^{(m+1)}(x,0)\\[2mm]
\tilde{a}^{(m+1)}(x) & \tilde{b}^{(m+1)}(x) & \tilde{m}_1^{(m+1)}(x)& \cdots & \tilde{m}_n^{(m+1)}(x) &\tilde{\ell}^{(m+1)}(x) \\[2mm]
R^{(m+2)}(x,0) & R_t^{(m+2)}(x,0) &\dis\pa_{x_1} R^{(m+2)}(x,0) & \cdots & \dis\pa_{x_n} R^{(m+2)}(x,0) & \Delta R^{(m+2)}(x,0)\\[2mm]
\end{bmatrix} 
\end{equation}
with $\tilde{a}^{(i)}$, $\tilde{b}^{(i)}$, $\tilde{m}_k^{(i)}$ and $\tilde{\ell}^{(i)}$ defined as in (\ref{rndcoef}). Again since elementary row operations do not change the determinant, $\widetilde{U}_{q_1}(x)$ will have the same determinant as the matrix $\widetilde{U}(x)$ in the assumption (\ref{positivity2}). This completes the proof of Theorem \ref{Th3}.

\medskip
\noindent{\bf Proof of Theorem \ref{Th4}}. After achieving the uniqueness for the inverse source problem, we now prove the corresponding stability estimate (\ref{stability}). The proof below works essentially for both of the cases whether $n$ is odd or even, the only difference is in the choices of the functions $R^{(i)}$, $i=1, \cdots, m+2$, as indicated in the Theorem \ref{Th3}. 
First we go back to the inequality (\ref{key}), integrate over $\Omega$ gives
\begin{equation}\label{intkey}
\|f_0\|_{L^2(\Omega)}^2+\|f_1\|_{L^2(\Omega)}^2+\|f_2\|_{L^2(\Omega)}^2+\|{\f}\|_{{\bf L}^2(\Omega)}^2 \leq C \sum_{i=1}^{m+2}\left(\|u^{(i)}_{tt}(\cdot,0)\|^2_{L^2(\Omega)}+\|u^{(i)}_{ttt}(\cdot,0)\|^2_{L^2(\Omega)}\right).
\end{equation}
 
For each $i$, $1\leq i\leq m+2$, we return to the $u_{tt}^{(i)}$-system:
\begin{equation}\label{th2utt}
\begin{cases}
(u^{(i)}_{tt})_{tt} - c^2(x) \Delta u^{(i)}_{tt}+q_1(x)(u^{(i)}_{tt})_t+q_0(x)u^{(i)}_{tt}+ {\q}(x)\cdot\nabla u^{(i)}_{tt} = S^{(i)}_{tt}(x,t)\\[2mm]
u^{(i)}_{tt}(x, 0)  = S^{(i)}(x,0), \ (u^{(i)}_{tt})_t(x, 0) = S^{(i)}_t(x,0)-q_1(x)S^{(i)}(x,0)  \\[2mm]
u^{(i)}_{tt}|_{\Gamma\times [-T, T]}=0 
\end{cases}
\end{equation}
with $S^{(i)}(x,t)=f_0(x)R^{(i)}+f_1(x)R_t^{(i)}+{\f}(x)\cdot \nabla R^{(i)}+f_2(x)\Delta R^{(i)}$. 
Here we assume 
\be\label{54}
c\in\mathcal{C}, \ q_0, q_1, q_2\in L^{\infty}(\Omega), \ {\q}\in \left(L^{\infty}(\Omega)\right)^n, f_0, f_1, f_2\in H^1_0(\Omega), \ {\f}\in \left(H^{1}_0(\Omega)\right)^n
\ee
 and $R^{(i)}$ satisfies (\ref{regularity2}) and (\ref{positivity1}) (or (\ref{positivity2}) if $n$ is even).
By linearity, we split $u_{tt}^{(i)}$ into two systems, $u_{tt}^{(i)} = y^{(i)}+z^{(i)}$, where $y^{(i)}=y^{(i)}(x,t)$ satisfies the homogeneous forcing term and nonhomogeneous initial conditions
\begin{equation}\label{th2y}
\begin{cases}
y^{(i)}_{tt} - c^2(x)\Delta y^{(i)}+ q_1(x)y^{(i)}_t+q_0(x)y^{(i)}+ {\q}(x)\cdot\nabla y^{(i)} = 0 & \mbox{in } \ Q \\[2mm]
y^{(i)}(x, 0) = u^{(i)}_{tt}(x, 0)  = S^{(i)}(x,0) & \mbox{in } \ \Omega\\[2mm]
y_t^{(i)}(x,0) =(u^{(i)}_{tt})_t(x, 0) = S^{(i)}_t(x,0)-q_1(x)S^{(i)}(x,0) & \mbox{in } \ \Omega \\[2mm]
y^{(i)}|_{\Gamma\times [-T, T]}=0 & \mbox{in } \ \Sigma
\end{cases}
\end{equation}
and $z^{(i)}=z^{(i)}(x,t)$ has the nonhomogeneous forcing term and homogeneous initial conditions
\begin{equation}\label{th2z}
\begin{cases}
z^{(i)}_{tt} - c^2(x)\Delta z^{(i)}+ q_1(x)z^{(i)}_t+q_0(x)z^{(i)}+ {\q}(x)\cdot\nabla z^{(i)} = S^{(i)}_{tt}(x,t) & \mbox{in } \ Q \\[2mm]
z^{(i)}(x, 0) =z^{(i)}_t(x, 0) = 0  & \mbox{in } \ \Omega \\[2mm]
z^{(i)}|_{\Gamma\times [-T, T]}=0 & \mbox{in } \ \Sigma.
\end{cases}
\end{equation}

For the $y^{(i)}$-system, note by assumptions (\ref{54}) and (\ref{regularity2}) we have 
\begin{equation*}
S^{(i)}(\cdot,0)\in H^1_0(\Omega) \ \ \mbox{and} \ \ S^{(i)}_t(\cdot,0)-q_1(\cdot)S^{(i)}(\cdot,0)\in L^2(\Omega). 
\end{equation*}
Thus we may apply the continuous observability inequality (\ref{coi}) (with $g=c^{-2}(x)dx^2$) to get 
\begin{equation*}
\|y^{(i)}(\cdot,0)\|_{H_0^{1}(\Omega)}^2 + \|y_t^{(i)}(\cdot,0)\|_{L^2(\Omega)}^2=\|u^{(i)}_{tt}(\cdot,0)\|_{H_0^{1}(\Omega)}^2 + \|u^{(i)}_{ttt}(\cdot,0)\|_{L^2(\Omega)}^2 \leq C \left\|\frac{\pa y^{(i)}}{\pa\nu}\right\|^2_{L^2{(\Sigma_1)}}.
\end{equation*}

Sum the above inequality over $i$, use (\ref{intkey}) and the decomposition $u_{tt}^{(i)} = y^{(i)}+z^{(i)}$, as well as Poincar\'{e}'s inequality, we have
\begin{eqnarray}\label{drop}
& \ & \|f_0\|_{L^2(\Omega)}^2+\|f_1\|_{L^2(\Omega)}^2+\|f_2\|_{L^2(\Omega)}^2+\|{\f}\|_{{\bf L}^2(\Omega)}^2 \\[2mm] 
& \leq & C \sum_{i=1}^{m+2}\left(\|u^{(i)}_{ttt}(\cdot,0)\|_{H^1_0(\Omega)}^2+\|u^{(i)}_{ttt}(\cdot,0)\|_{L^2(\Omega)}^2\right)\nonumber\\[2mm]
& \leq & C \sum_{i=1}^{m+2}\left\|\frac{\pa y^{(i)}}{\pa\nu}\right\|^2_{L^2{(\Sigma_1)}} \nonumber\\[2mm] 
& = & C \sum_{i=1}^{m+2}\left\|\frac{\pa u_{tt}^{(i)}}{\pa\nu} - \frac{\pa z^{(i)}}{\pa\nu}\right\|^2_{L^2{(\Sigma_1)}} \nonumber \\[2mm]
& \leq & C \sum_{i=1}^{m+2}\left\|\frac{\pa u_{tt}^{(i)}}{\pa\nu}\right\|^2_{L^2{(\Sigma_1)}} + C \sum_{i=1}^{m+2}\left\|\frac{\pa z^{(i)}}{\pa\nu}\right\|^2_{L^2{(\Sigma_1)}}.\nonumber
\end{eqnarray}

Note this is the desired stability estimate (\ref{stability}) polluted by the $z^{(i)}$ terms. Next we show those terms can be absorbed through a compactness--uniqueness argument, where the uniqueness relies on Theorem \ref{Th3}. To start, note for the $z^{(i)}$-system (\ref{th2z}), we have the following proposition.
\begin{proposition}
For each $i=1, \cdots, m+2$, the operator define by 
\begin{eqnarray}
\calK_i: L^2(\Omega)\times L^2(\Omega)\times L^2(\Omega)\times\left(L^2(\Omega)\right)^n&\to&\ L^2(\Sigma_1)  \\[2mm]
 (f_0, f_1, f_2, {\f }) &\mapsto &\frac{\pa z^{(i)}}{\pa\nu}|_{\Sigma_1}, \nonumber
\end{eqnarray}
is a compact operator.
\end{proposition}  
\begin{proof}
Note assumptions (\ref{54}) and (\ref{regularity2}) imply $S^{(i)}_{tt}\in H^{1}(Q)$, thus by the regularity result (\ref{reg}) we have
\begin{equation*}
S^{(i)}_{tt}\in H^{1}(Q) \ \Rightarrow \ \frac{\pa z^{(i)}}{\pa\nu}\in H^{1}(\Sigma_1) \ \ \mbox{continuously}.
\end{equation*}
This then implies the map $(f_0, f_1, f_2, {\f })\mapsto \calK_{i}(f_0, f_1, f_2, {\f })\in H^1(\Sigma_1)$ is continuous and hence $(f_0, f_1, f_2, {\f })\mapsto \calK_{i}(f_0, f_1, f_2, {\f })\in L^2(\Sigma_1)$ is compact.
\end{proof}

Combine $\calK^{(i)}$, $i=1, \cdots, m+2$, being compact, together with the uniqueness result in Theorem \ref{Th3}, we may drop the $z^{(i)}$ terms in (\ref{drop}) to get the desired stability estimate (\ref{stability}). To carry this out, suppose by contradiction the stability estimate (\ref{stability}) does not hold, then there exist sequences $\{f_0^k\}$, $\{f_1^k\}$, $\{f_2^k\}$ and $\{{\f}^k\}$, with $f_0^k, f_1^k, f_2^k\in H^1_0(\Omega)$ and ${\f}^k\in \left(H^1_0(\Omega)\right)^n$, $\forall k\in\mathbb{N}$, such that 
\begin{equation}\label{hyp1}
\left\|f_0^k\right\|_{L^2(\Omega)}^2+\left\|f_1^k\right\|_{L^2(\Omega)}^2+\left\|f_2^k\right\|_{L^2(\Omega)}^2+\left\|{\f}^k\right\|_{{\bf L}^2(\Omega)}^2 = 1
\end{equation}
and
\begin{equation}\label{hyp2}
\lim_{k\to\infty}\sum_{i=1}^{m+2}\left\|\frac{\pa u_{tt}^{(i)}(f_0^k, f_1^k, f_2^k,{\f}^k)}{\pa\nu}\right\|_{L^2(\Sigma_1)}=0
\end{equation}
where $u^{(i)}(f_0^k, f_1^k, f_2^k,{\f}^k)$ solves the system (\ref{th1u}) with $f_0=f_0^k$, $f_1=f_1^k$, $f_2=f_2^k$ and ${\f} = {\f}^k$. 
From (\ref{hyp1}), there exist subsequences, still denoted as $\{f_0^k\}$, $\{f_1^k\}$, $\{f_2^k\}$ and $\{{\f}^k\}$, such that 
\begin{equation}\label{weak}
{f_i}^k \rightharpoonup f_i^* \ \mbox{and} \ {\f}^k \rightharpoonup {\f}^* \ \mbox{weakly for some} \ f_i^*\in L^2(\Omega) \ \mbox{and} \ {\f}^*\in \left(L^2(\Omega)\right)^n, i=0,1,2.
\end{equation}
Moreover, in view of the compactness of $\calK_{i}$, $i=1,\cdots, m+2$, we also have the strong convergence
\begin{equation}\label{strong}
\lim_{k, l\to\infty}\left\|\calK_i(f_0^k, f_1^k, f_2^k,{\f}^k) - \calK_i(f_0^l, f_1^l, f_2^l,{\f}^l)\right\|_{L^2(\Sigma_1)} = 0, \ \forall i=1, \cdots, m+2.
\end{equation}

On the other hand, since the map $(f_0,f_1,f_2,{\f})\mapsto u^{(i)}(f_0,f_1,f_2,{\f})$ is linear, we have from (\ref{drop}) that
\begin{eqnarray*}
& \ & \left\|f_0^k-f_0^l\right\|^2_{L^2(\Omega)}+\left\|f_1^k-f_1^l\right\|^2_{L^2(\Omega)}+\left\|f_2^k-f_2^l\right\|^2_{L^2(\Omega)}+\left\|{\f}^k-{\f}^l\right\|^2_{{\bf L}^2(\Omega)} \\[2mm]
& \leq & C \sum_{i=1}^{m+2}\left\|\frac{\pa u_{tt}^{(i)}(f_0^k, f_1^k, f_2^k,{\f}^k)}{\pa\nu}-\frac{\pa u_{tt}^{(i)}(f_0^l, f_1^l, f_2^l,{\f}^l)}{\pa\nu}\right\|^2_{L^2(\Sigma_1)} \\
& \ & +C\sum_{i=1}^{m+2}\left\|\calK_i(f_0^k, f_1^k, f_2^k,{\f}^k) - \calK_i(f_0^l, f_1^l, f_2^l,{\f}^l)\right\|^2_{L^2(\Sigma_1)} \nonumber\\[2mm]
& \leq & C \sum_{i=1}^{m+2}\left\|\frac{\pa u_{tt}^{(i)}(f_0^k, f_1^k, f_2^k,{\f}^k)}{\pa\nu}\right\|^2_{L^2(\Sigma_1)}+C \sum_{i=1}^{m+2}\left\|\frac{\pa u_{tt}^{(i)}(f_0^l, f_1^l, f_2^l,{\f}^l)}{\pa\nu}\right\|^2_{L^2(\Sigma_1)} \\
& \ &  +C\sum_{i=1}^{m+2}\left\|\calK_i(f_0^k, f_1^k, f_2^k,{\f}^k) - \calK_i(f_0^l, f_1^l, f_2^l,{\f}^l)\right\|^2_{L^2(\Sigma_1)} \nonumber
\end{eqnarray*}
and therefore by (\ref{hyp2}) and (\ref{strong}) we get 
\begin{equation*}
\lim_{k, l\to\infty}\left\|f_i^k-f_i^l\right\|_{L^2(\Omega)}=\lim_{k, l\to\infty}\left\|{\f}^k-{\f}^l\right\|_{{\bf L}^2(\Omega)} = 0, \ i=0,1,2.
\end{equation*}
Namely, $\{f_0^k\}$, $\{f_1^k\}$, $\{f_2^k\}$ are Cauchy sequences in $L^2(\Omega)$ and $\{{\f}_k\}$ is a Cauchy sequence in $(L^2(\Omega))^n$. By uniqueness of limit and in view of (\ref{weak}), we must have 
$\{f_i^k\}$ converges to $f_i^*$ strongly, $i=0,1,2$, and $\{{\f}^k\}$ converges to ${\f}^*$ strongly. Hence we have from (\ref{hyp1})
\be\label{63}
\|f_0^*\|_{L^2(\Omega)}^2+\|f_1^*\|_{L^2(\Omega)}^2+\|f_2^*\|_{L^2(\Omega)}^2+\|{\f}^*\|_{{\bf L}^2(\Omega)}^2 =1.
\ee

Now again for the $u_{tt}^{(i)}$-system (\ref{th2utt}), by the regularity theory (\ref{reg}) we have that the map $(f_0,f_1,f_2,{\f})\mapsto \frac{\pa u_{tt}^{(i)}(f_0,f_1,f_2,{\f})}{\pa\nu}\in L^2(\Sigma)$ is continuous and hence 
\begin{equation*}
\left\|\frac{\pa u_{tt}^{(i)}(f_0,f_1,f_2,{\f})}{\pa\nu}\right\|_{L^2(\Sigma)}^2\leq C \left(\|f_0\|_{L^2(\Omega)}^2+\|f_1\|_{L^2(\Omega)}^2+\|f_2\|_{L^2(\Omega)}^2+\|{\f}\|_{{\bf L}^2(\Omega)}^2\right).
\end{equation*}
Since the map $(f_0,f_1,f_2,{\f})\mapsto u_{tt}^{(i)}(f_0,f_1,f_2,{\f})|_{\Sigma}$ is linear, we thus have 
\begin{eqnarray}
& \ & \left\|\frac{\pa u_{tt}^{(i)}(f_0^k,f_1^k,f_2^k,{\f}^k)}{\pa\nu}-\frac{\pa u_{tt}^{(i)}(f_0^*,f_1^*,f_2^*,{\f}^*)}{\pa\nu}\right\|_{L^2(\Sigma_1)}^2 \\[2mm]
& \leq & C \left(\|f_0^k-f_0^*\|_{L^2(\Omega)}^2+\|f_1^k-f_1^*\|_{L^2(\Omega)}^2+\|f_2^k-f_2^*\|_{L^2(\Omega)}^2+\|{\f}^k-{\f}^*\|_{{\bf L}^2(\Omega)}^2\right). \nonumber
\end{eqnarray}
This then implies, by virtue of $f_i^k\to f_i^*$, $i=0,1,2$ and ${\bf f}^k\to {\bf f}^*$ strongly, that  
\begin{equation*}
\lim_{k\to\infty}\left\|\frac{\pa u_{tt}^{(i)}(f_0^k,f_1^k,f_2^k,{\f}^k)}{\pa\nu}-\frac{\pa u_{tt}^{(i)}(f_0^*,f_1^*,f_2^*,{\f}^*)}{\pa\nu}\right\|_{L^2(\Sigma_1)} = 0
\end{equation*}
and hence $\displaystyle\frac{\pa u_{tt}^{(i)}(f_0^*,f_1^*,f_2^*,{\f}^*)}{\pa\nu} = 0$ in $L^2(\Sigma_1)$ in view of (\ref{hyp2}). In other words, $\displaystyle\frac{\pa u_{t}^{(i)}(f_0^*,f_1^*,f_2^*,{\f}^*)}{\pa\nu}$ is a constant in $t\in[-T,T]$.
We claim that $\displaystyle\frac{\pa u_{t}^{(i)}(f_0^*,f_1^*,f_2^*,{\f}^*)}{\pa\nu}=0$ on $\Sigma_1$. To see this, we consider the $u_t^{(i)}(f_0^k,f_1^k,f_2^k,{\f}^k)$-system
\begin{equation}\label{th2ut}
\begin{cases}
(u^{(i)}_t)_{tt} - c^2(x)\Delta u^{(i)}_t +q_1(x)(u_t^{(i)})_t+q_0(x)u^{(i)}_t+{\q}(x)\cdot\nabla u^{(i)}_t=(S_k^{(i)})_t(x,t) & \mbox{in } \ Q \\[2mm]
(u_t^{(i)})(x, 0)=0, \ (u^{(i)}_t)_t(x, 0) = S_k^{(i)}(x,0) & \mbox{in } \ \Omega  \\[2mm]
u^{(i)}_{t}|_{\Gamma\times [-T, T]}=0  & \mbox{in } \ \Sigma
\end{cases}
\end{equation}
where for $i = 1, \cdots, m+2$, $R^{(i)}=R^{(i)}(x,t)$ and
\begin{equation*}
S_k^{(i)}(x,t) = f_0^k(x)R^{(i)}+f_1^k(x)R_t^{(i)}+{\f}^k(x)\cdot \nabla R^{(i)}+f_2^k(x)\Delta R^{(i)}.
\end{equation*}
The standard regularity theory (\ref{reg}) and trace theory implies
\begin{eqnarray*}
& \ & \left\|u_t^{(i)}(f_0^k,f_1^k,f_2^k,{\f}^k) - u_t^{(i)}(f_0^*,f_1^*,f_2^*,{\f}^*)\right\|_{C([0,T]; H_0^1(\Omega))}^2 \\[2mm]
& \leq & C \left(\|f_0^k-f_0^*\|_{L^2(\Omega)}^2+\|f_1^k-f_1^*\|_{L^2(\Omega)}^2+\|f_2^k-f_2^*\|_{L^2(\Omega)}^2+\|{\f}^k-{\f}^*\|_{{\bf L}^2(\Omega)}^2\right)
\end{eqnarray*}
and 
\begin{eqnarray*}
& \ & \left\|u_t^{(i)}(f_0^k,f_1^k,f_2^k,{\f}^k) - u_t^{(i)}(f_0^*,f_1^*,f_2^*,{\f}^*)\right\|_{C([0,T]; H^{\frac{1}{2}}(\Sigma)}^2 \\[2mm]
& \leq & C \left(\|f_0^k-f_0^*\|_{L^2(\Omega)}^2+\|f_1^k-f_1^*\|_{L^2(\Omega)}^2+\|f_2^k-f_2^*\|_{L^2(\Omega)}^2+\|{\f}^k-{\f}^*\|_{{\bf L}^2(\Omega)}^2\right).
\end{eqnarray*}

Note $u_t^{(i)}(f_0^k,f_1^k,f_2^k,{\f}^k)(x, 0)=0$ as well as the strong convergence of $f_i^k\to f_i^*$, $i=0,1,2$ and ${\bf f}^k\to {\bf f}^*$. Thus letting $k\to\infty$ we get $u_t^{(i)}(f_0^*,f_1^*,f_2^*,{\f}^*)(x, 0)=0$ in $\Omega$ and $u_t^{(i)}(f_0^*,f_1^*,f_2^*,{\f}^*)|_{\Sigma}=0$. Hence $\displaystyle\frac{\pa u_{t}^{(i)}(f_0^*,f_1^*,f_2^*,{\f}^*)}{\pa\nu}(x,0) = 0$ on $\Sigma$. Since we know $\displaystyle\frac{\pa u_{t}^{(i)}(f_0^*,f_1^*,f_2^*,{\f}^*)}{\pa\nu}$ is a constant in $t$, we must have $\displaystyle\frac{\pa u_{t}^{(i)}(f_0^*,f_1^*,f_2^*,{\f}^*)}{\pa\nu}=0$ on $\Sigma_1$, as desired.

The above then implies $\displaystyle\frac{\pa u^{(i)}(f_0^*,f_1^*,f_2^*,{\f}^*)}{\pa\nu}$ is also a constant in $t$. By repeating the same argument, this time using the regularity theory for the $u^{(i)}(f_0^k,f_1^k,f_2^k,{\f}^k)$-system and taking limit $k\to\infty$, we finally get $\displaystyle\frac{\pa u^{(i)}(f_0^*,f_1^*,f_2^*,{\f}^*)}{\pa\nu}=0$ on $\Sigma_1$.

Hence we have that $u^{(i)}(f_0^*,f_1^*,f_2^*,{\f}^*)$ satisfies the following 
\begin{equation}\label{th2uf0}
\begin{cases}
u^{(i)}_{tt} - c^2(x)\Delta u^{(i)} + q_1(x)u_t^{(i)}+q_0(x)u^{(i)}+{\q}(x)\cdot\nabla u^{(i)} = S_*^{(i)}(x,t) & \mbox{in } \ Q \\[2mm]
u^{(i)}(x, 0)  = u^{(i)}_t(x, 0) = 0 & \mbox{in } \ \Omega  \\[2mm]
u^{(i)}|_{\Gamma\times[-T,T]}=0, \quad \frac{\pa u^{(i)}}{\pa\nu}|_{\Gamma_1\times[-T,T]}=0 & \mbox{in } \ \Sigma, \Sigma_1
\end{cases}
\end{equation}
with $$S_*^{(i)}(x,t)=f_0^*(x)R^{(i)}+f_1^*(x)R_t^{(i)}+{\f}^*(x)\cdot \nabla R^{(i)}+f_2^*(x)\Delta R^{(i)}, i = 1, \cdots, m+2.$$
By the uniqueness result we proved in Theorem {\ref{Th3}, this must imply $f_0^*=f_1^*=f_2^*={\f}^*= 0$, which contradicts with (\ref{63}). 
Hence we must be able to drop the $z^{(i)}$ terms in (\ref{drop}). This completes the proof of Theorem \ref{Th4}.

\medskip
{\bf Proof of Theorems \ref{Th1} and \ref{Th2}}. Finally, we provide the proofs of uniqueness and stability of the original inverse problem. These results are pretty much direct consequences of Theorems \ref{Th3} and \ref{Th4} given the relationship (\ref{relation}) between the original inverse problem and the inverse source problem. More precisely, we have the positivity conditions (\ref{positivity1'}) and (\ref{positivity2'}) imply (\ref{positivity1}) and (\ref{positivity2}). In addition, by the regularity theory (\ref{reg}) the assumption (\ref{regularity}) on the initial and boundary conditions $\{w_0^{(i)}, w_1^{(i)}, h\}$ implies the solutions $w^{(i)}$, $i=1, \cdots, m+2$, satisfy 
\begin{equation*}
\{w^{(i)}, w_t^{(i)}, w_{tt}^{(i)}, w_{ttt}^{(i)}\}\in C\left([-T,T]; H^{\gamma+1}(\Omega)\times H^{\gamma}(\Omega)\times H^{\gamma-1}(\Omega)\times H^{\gamma-2}(\Omega)\right).
\end{equation*} 
As $\gamma>\frac{n}{2}+4$, we have the following embedding $H^{\gamma-2}(\Omega)\hookrightarrow W^{2,\infty}(\Omega)$ and hence the regularity assumption (\ref{regularity}) implies the corresponding regularity assumption (\ref{regularity2}) for the inverse source problem. This completes the proof of all the theorems.

\section{Some Examples and Concluding Remarks}

In this last section we first provide some concrete examples such that the key positivity conditions (\ref{positivity1'}), (\ref{positivity2'}), (\ref{positivity1}), (\ref{positivity2}) are satisfied, and then give some general remarks.
 
\medskip
\noindent{\bf Example 1}. Consider the following functions $R^{(i)}(x,t)$, $x=(x_1, \cdots, x_n)\in\Omega$, $t\in [-T, T]$, $i=1, \cdots, m+2$, defined by
\begin{eqnarray*}
R^{(1)}(x,t) & = &  t, \ R^{(i)}(x,t) = x_{2i-3}+tx_{2i-2}, \ 2\leq i\leq m+1, \\[2mm]
R^{(m+2)}(x,t) & = & \begin{cases} x_{2m+1}+\frac{1}{2}tx_1^2 & \ \mbox{if} \ n=2m+1 \ \mbox{is odd} \\[2mm] \frac{1}{2}x_1^2 & \ \mbox{if} \ n=2m \ \mbox{is even}.\end{cases}
\end{eqnarray*} 
Then we may easily see that the matrices $U(x)$ and $\widetilde{U}(x)$ are lower triangular matrices with all $1$s at the diagonal after swapping the first two columns. Thus the determinants of the matrices $U(x)$ and $\widetilde{U}(x)$ are both $-1$ and hence conditions (\ref{positivity1}), (\ref{positivity2}) are satisfied. Correspondingly, we may choose the $m+2$ pairs of initial conditions $\{w_0^{(i)}(x), w_1^{(i)}(x)\}$ as 
\begin{eqnarray*}
w_0^{(1)}(x) & = & 0,   \ w_1^{(1)}(x)=1, \\[2mm]  
w_0^{(i)}(x) & = & x_{2i-3}, \ w_1^{(i)}(x)=x_{2i-2}, \ 2\leq i\leq m+1, \\[2mm]
w_0^{(m+2)}(x) & = & \begin{cases} x_{2m+1} & \ \mbox{if} \ n=2m+1 \ \mbox{is odd} \\[2mm] \frac{1}{2}x_1^2 & \ \mbox{if} \ n=2m \ \mbox{is even}\end{cases} \\[2mm]
w_1^{(m+2)}(x) & = & \begin{cases} \frac{1}{2}x_1^2 & \ \mbox{if} \ n=2m+1 \ \mbox{is odd} \\[2mm] 0 & \ \mbox{if} \ n=2m \ \mbox{is even}.\end{cases}
\end{eqnarray*} 
Then the matrices $W(x)$ and $\widetilde{W}(x)$ are also lower triangular matrices with all $1$s at the diagonal after swapping the first two columns and hence conditions are (\ref{positivity1'}), (\ref{positivity2'}) are satisfied.

\medskip
\noindent{\bf Example 2}. Considering the following functions $R^{(i)}(x,t)$, $x=(x_1, \cdots, x_n)\in\Omega$, $t\in [-T, T]$, $i=1, \cdots, m+2$, defined by
\begin{eqnarray*}
R^{(1)}(x,t) & = &  \sin{t}, \ R^{(i)}(x,t) = \cos{t}e^{x_{2i-3}}+\sin{t}e^{x_{2i-2}}, \ 2\leq i\leq m+1, \\[2mm]
R^{(m+2)}(x,t) & = & \begin{cases} \cos{t}e^{x_{2m+1}}+\sin{t}e^{-x_1} & \ \mbox{if} \ n=2m+1 \ \mbox{is odd} \\[2mm] \cos{t}e^{-x_1} & \ \mbox{if} \ n=2m \ \mbox{is even}.\end{cases}
\end{eqnarray*} 

Then the matrix $U(x)$ becomes
\begin{align}
    U(x)=\begin{bmatrix}
        0 & 1 & 0 & 0 & \cdots & 0 & 0 & 0 \\[2mm]
        1 & 0 & 0 & 0 & \cdots & 0 & 0 & 0 \\[2mm]
        e^{x_1} & e^{x_2} & e^{x_1} & 0 & 0 & \cdots & 0 & e^{x_1} \\[2mm]
        e^{x_2} &-e^{x_1} & 0 & e^{x_2} & 0 & \cdots & 0 & e^{x_2}\\[2mm]
        \vdots & \vdots & \vdots & \ddots & \ddots & \ddots & \vdots & \vdots \\
        e^{x_{n-1}} &-e^{x_{n-2}} & 0 & \cdots & 0 & e^{x_{n-1}} &0& e^{x_{n-1}}\\[2mm]
        e^{x_{n}} & e^{-x_1} & 0 & \cdots & 0 & 0 & e^{x_n}  & e^{x_n} \\[2mm]
        e^{-x_1} &-e^{x_n} &-e^{-x_1} & 0 & \cdots & 0 & 0 & e^{-x_1}\\[2mm]
    \end{bmatrix} 
\end{align}
Notice that $U(x)$ is not a lower triangular matrix. However, we can easily transform the it into a lower triangular matrix by swapping the first two columns and subtracting the 3rd, 4th, ..., (n+2)th column from the last column. As a consequence we get $\det{U}(x) = -2\prod_{i=2}^ne^{x_i}$. In a similar fashion we can also get $\det{\tilde{U}}(x) = -2\prod_{i=2}^ne^{x_i}$. As $\Omega$ is a bounded domain, we hence have the conditions (\ref{positivity1}) and (\ref{positivity2}) are satisfied.

Correspondingly we may choose the $m+2$ pairs of initial conditions $\{w_0^{(i)},w_1^{(i)}\}$ as 
\begin{eqnarray*}
w_0^{(1)}(x) & = & 0,   \ w_1^{(1)}(x)=1, \\[2mm]  
w_0^{(i)}(x) & = & e^{x_{2i-3}}, \ w_1^{(i)}(x)=e^{x_{2i-2}}, \ 2\leq i\leq m+1, \\[2mm]
w_0^{(m+2)}(x) & = & \begin{cases}e^{x_{n}} & \ \mbox{if} \ n=2m+1 \ \mbox{is odd} \\[2mm] e^{-x_1} & \ \mbox{if} \ n=2m \ \mbox{is even}\end{cases} \\[2mm]
w_1^{(m+2)}(x) & = & \begin{cases} e^{-x_1} & \ \mbox{if} \ n=2m+1 \ \mbox{is odd} \\[2mm] 0 & \ \mbox{if} \ n=2m \ \mbox{is even}.\end{cases}
\end{eqnarray*} 
Then the determinants of both $W(x)$ and $\widetilde{W}(x)$ are also $-2\prod_{i=2}^ne^{x_i}$, calculated in the same manner as in the case of $U(x)$ and $\widetilde{U}(x)$. Hence the conditions (\ref{positivity1'}) and (\ref{positivity2'}) are satisfied.

\medskip
\noindent{\bf Example 3}.
In general if we have $f^{(j)}\in C^2(\overline{\Omega})$ with $f^{(j)}(x)=f^{(j)}(x_1, \cdots, x_j)$, $1\leq j\leq n$ and $g, h\in C^2[-T,T]$ that satisfies
\begin{eqnarray*}
\left|\frac{\partial f^{(j)}}{\partial x_j}\right|\geq r_j>0, 1\leq j\leq n, \qquad \left|\frac{\partial^2 f^{(1)}}{\partial x_1^2}\right|\geq \tilde{r}_1>0 \\[2mm]
g(0)=h'(0)=1, \quad g'(0)=h(0)=0.
\end{eqnarray*}
for some positive $r_j$, $1\leq j\leq n$, and $\tilde{r}_1$. Then we may consider the functions $R^{(i)}(x,t)$, $x=(x_1, \cdots, x_n)\in\Omega$, $t\in [-T, T]$, $i=1, \cdots, m+2$, of the following form:
\begin{eqnarray*}
R^{(1)}(x,t) & = & h(t), \ R^{(i)}(x,t) = f^{(2i-3)}(x)g(t)+f^{(2i-2)}(x)h(t), \ 2\leq i\leq m+1, \\[2mm]
R^{(m+2)}(x,t) & = & \begin{cases} f^{(n)}(x)g(t)+f^{(1)}(ax)h(t) & \ \mbox{if} \ n=2m+1 \ \mbox{is odd} \\[2mm] f^{(1)}(ax)g(t) & \ \mbox{if} \ n=2m \ \mbox{is even}.\end{cases}
\end{eqnarray*} 
where $a<0$ so that $ax\in\Omega$. 

Correspondingly we may choose the $m+2$ pairs of initial conditions $\{w_0^{(i)},w_1^{(i)}\}$ as 
\begin{eqnarray*}
w_0^{(1)}(x) & = & 0,   \ w_1^{(1)}(x)=1, \\[2mm]  
w_0^{(i)}(x) & = & f^{(2i-3)}(x), \ w_1^{(i)}(x)=f^{(2i-2)}(x), \ 2\leq i\leq m+1, \\[2mm]
w_0^{(m+2)}(x) & = & \begin{cases} f^{(n)}(x) & \ \mbox{if} \ n=2m+1 \ \mbox{is odd} \\[2mm] f^{(1)}(ax) & \ \mbox{if} \ n=2m \ \mbox{is even}\end{cases} \\[2mm]
w_1^{(m+2)}(x) & = & \begin{cases} f^{(1)}(ax) & \ \mbox{if} \ n=2m+1 \ \mbox{is odd} \\[2mm] 0 & \ \mbox{if} \ n=2m \ \mbox{is even}.\end{cases}
\end{eqnarray*} 

In this case, after swapping the first and second column, the last column with the preceding $(n+2)$th, $(n+1)$th, $\cdots$, and finally the $3$rd column, as well as swapping the last row with the preceding $(n+2)$th, $(n+1)$th, $\cdots$, and finally the $3$rd row. We may get the determinants of the matrices 
$U(x)$, $\widetilde{U}(x)$, $W(x)$ and $\widetilde{W}(x)$ are equal to 
$$ \left(a\partial_{x_1}f^{(1)}({a}x)\partial_{x_1}^2f^{(1)}(x)-{a}^2\partial_{x_1}^2f^{(1)}({a}x)\partial_{x_1}f^{(1)}(x)\right)\prod^n_{j=2}\partial_{x_j}f^{(j)}(x).$$

Since $f^{(1)}\in C^2(\overline{\Omega})$, $|\partial_{x_1}f^{(1)}|\geq r_1>0$ and $|\partial^2_{x_1}f^{(1)}|\geq \tilde{r}_1>0$, $\partial_{x_1}f^{(1)}$ and $\partial^2_{x_1}f^{(1)}$ do not change sign. Hence we have 
\begin{equation*}
|\left(a\partial_{x_1}f^{(1)}({a}x)\partial_{x_1}^2f^{(1)}(x)-{a}^2\partial_{x_1}^2f^{(1)}({a}x)\partial_{x_1}f^{(1)}(x)\right)\prod^n_{j=2}\partial_{x_j}f^{(j)}(x)|\geq (a^2+|a|)\tilde{r}_1\prod^n_{j=1}r_j.    
\end{equation*}
Hence the positivity conditions (\ref{positivity1'}), (\ref{positivity2'}), (\ref{positivity1}) and (\ref{positivity2}) are satisfied.

\bigskip
Finally, we end the paper with some comments and remarks.

\medskip
(1) We have shown in this paper that in order to recover all the coefficients, we need to appropriately choose $\lfloor{\frac{n+4}{2}}\rfloor$ pairs of initial conditions $\{w_0, w_1\}$ and a boundary condition $h$, and then use their corresponding boundary measurements. As mentioned earlier, since in total there are  $n+3$ unknown functions, it is natural to expect to recover them from $n+3$ boundary measurements. Indeed, following the approach of this paper, we can also achieve the recovery by appropriately choosing $n+3$ initial positions $w_0$ with an initial velocity $w_1$ and a boundary condition $h$, and then use their corresponding boundary  measurements. In particular, in this case the positivity condition becomes
\begin{equation*}
\det 
\begin{bmatrix}
w_0^{(1)}(x) & w_1(x) & \dis\pa_{x_1} w_0^{(1)}(x) & \cdots & \dis\pa_{x_n} w_0^{(1)}(x) &\Delta w_0^{(1)}(x) \\[2mm]
w_0^{(2)}(x) & w_1(x) & \dis\pa_{x_1} w_0^{(2)}(x) & \cdots & \dis\pa_{x_n} w_0^{(2)}(x) &\Delta w_0^{(2)}(x) \\[2mm]
\vdots & \vdots  & \vdots & \ddots & \vdots & \vdots \\
w_0^{(n+3)}(x) & w_1(x) &\dis\pa_{x_1} w_0^{(n+3)}(x) & \cdots & \dis\pa_{x_n} w_0^{(n+3)}(x) & \Delta w_0^{(n+3)}(x)
\end{bmatrix}
\geq r_0 > 0.
\end{equation*}
Note although in this case we need more measurements, an advantage is that we only need to differentiate the $u$-equation with respect to $t$ twice, rather than three times. We may also get a better stability estimate of the form
\begin{multline*}
\|c^2-{\Tilde{c}}^2\|_{L^2(\Omega)}^2 + \|q_1-p_1\|_{L^2(\Omega)}^2+\|q_0-p_0\|_{L^2(\Omega)}^2+\|{\q}-{\p}\|^2_{{\bf L}^2(\Omega)}\\[2mm]
\leq C\sum_{i=1}^{n+3}\left\|\frac{\pa w_{t}^{(i)}({c,q_1,q_0,\q})}{\pa\nu}-\frac{\pa w_{t}^{(i)}({\tilde{c},p_1,p_0,\p})}{\pa\nu}\right\|_{L^2(\Sigma_1)}^2
\end{multline*}
in the sense that we only need to differentiate the measurements in time once.

\medskip
(2) In our problem formulation we use the time interval $[-T,T]$ and regard the middle 
$t=0$ as initial time. This is not essential since a simple change of variable $t\to t-T$ transforms
$t=0$ to $t=-T$. However, this present choice allows the recovery of all coefficients with fewer choices of initial conditions and hence fewer boundary measurements. This is because we may use both equations in (\ref{th1eq1}), compare to just one if we assume the time interval as $[0,T]$ and then extend solutions to $[-T,0]$.

\medskip
(3) It is also possible to set up the inverse problem by assuming Neumann boundary condition $\frac{\partial w}{\partial \nu}$ on $\Sigma=\Gamma\times [-T,T]$ and making measurements of Dirichlet boundary traces of the solution $w$ over $\Sigma_1=\Gamma_1\times[-T,T]$, such as in \cite{LiuT2011-2}. This, however, would require more demanding geometrical assumption on the unobserved portion of the boundary $\Gamma_0$. For example, we may need to assume $\frac{\pa d}{\pa\nu}=\langle Dd,\nu\rangle=0$ on $\Gamma_0$ in the geometrical assumption to account for the Neumann boundary condition \cite{TY2002}. In addition, the more delicate regularity theory of second-order hyperbolic equation with nonhomogeneous Neumann boundary condition will also need to be invoked \cite{LT1990}, \cite{LT1991}. Nevertheless, the main ideas of solving the inverse problem remain the same.

\bigskip
{\bf Acknowledgements}

\medskip
The first author would like to thank Professor Yang Yang for many very useful discussions.

\end{document}